\documentclass{commat}

\usepackage{graphicx}
\usepackage{multicol, multirow}

\newcommand{\A}{\mathcal{A}}
\newcommand{\C}{\mathbb{C}}
\newcommand{\K}{\mathbb{K}}
\newcommand{\N}{\mathbb{N}}
\newcommand{\Z}{\mathbb{Z}}
\DeclareMathOperator{\Der}{Der}

\title{%
On classification and deformations of Lie-Rinehart superalgebras
}

\author{%
Quentin Ehret and Abdenacer Makhlouf
}

\affiliation{
\address{Quentin Ehret --
Universit\'{e} de Haute-Alsace, IRIMAS UR 7499, F-68100 Mulhouse, France
}
\email{%
quentin.ehret@uha.fr
}
\address{Abdenacer Makhlouf --
Universit\'{e} de Haute-Alsace, IRIMAS UR 7499, F-68100 Mulhouse, France
}
\email{%
abdenacer.makhlouf@uha.fr
}
}

\abstract{%
The purpose of this paper is to study Lie-Rinehart superalgebras over characteristic zero fields, which are consisting of a supercommutative associative superalgebra $A$ and a Lie superalgebra $L$ that are compatible in a certain way. We discuss their structure and provide a classification in small dimensions.  We describe all  possible pairs defining a Lie-Rinehart superalgebra for $\dim(A)\leq 2$ and $\dim(L)\leq 4$. Moreover, we construct a cohomology complex and develop a theory of formal deformations based on formal power series and this cohomology. 
}



\VOLUME{30}
\NUMBER{2}
\firstpage{67}
\DOI{https://doi.org/10.46298/cm.10537}

\begin{paper}

\section*{Introduction}

Lie-Rinehart algebras are algebraic analogs of Lie algebroids. They first appeared in the work of Rinehart (\cite{RG63}) and Palais (\cite{PR61}) and have been studied by Huebschmann (\cite{HJ90} and \cite{HJ98}). A Lie-Rinehart algebra is a pair $(A,L)$, with $A$ an associative $\K$-algebra, $\K$ being a commutative ring, and $L$ a Lie $\K$-algebra. They must be endowed with an action of $A$ on $L$, making the latter an $A$-module, and with a Lie algebra map $\rho: L\longrightarrow\Der(A)$ such that $L$ acts on $A$ by derivations. Some authors put the emphasis on $L$ by saying that $L$ is a $\K$-$A$ Lie-Rinehart algebra. An example comes from Differential Geometry (\cite{RC20}): if $V$ is a differential manifold, let $A:=\mathcal{O}_V$ be the algebra of smooth functions on $V$, and $L:=\text{Vect}(V)$ be the Lie algebra of vector fields on $V$. Then the pair $(A,L)$ carries a Lie-Rinehart algebra structure. According to Claude Roger, Lie-Rinehart algebras give direct methods to deal with the relations between differential operators on manifolds and the enveloping algebra of the Lie algebra of vector fields of a manifold (\cite{RC20}). More details can be found in (\cite{RG63}).

The notion of Lie-Rinehart superalgebra, which is a generalization of Lie-Rinehart algebras to the superworld, has been studied for example by Chemla (\cite{CS95}) and Roger (\cite{RC20}). This theory can be used to study superalgebras of differential operators on supermanifolds, to generalize enveloping superalgebras related to those supermanifolds, and can have some applications in supersymmetry. Some examples can be found again in \cite{RC20}.

In this paper, we aim to explore and understand the $\C$-Lie-Rinehart superstructures in low dimensions. We give a classification of all the possible structures on $(A,L)$ when $\dim(A)\leq2$ and $\dim(L)\leq 4$. We rely heavily on the classification of associative and Lie superalgebras already existing in the literature (\cite{AM14}, \cite{ACZ09} and \cite{GP74} for the associative case; \cite{BN78} and \cite{BS99} for the Lie case) and on the computer algebra system Mathematica. We limit ourselves to those low dimensions due to technical and time constraints, even our program should work in any dimension. However, it appears that the number of possible structures are growing exponentially as soon as we increase the dimension, so it does not seem reasonable to give a list of all the possible Lie-Rinehart superalgebras in higher dimensions.

Another aspect we pursue in this paper is to provide a deformation theory and study one-parameter formal deformations of Lie-Rinehart superalgebras. In a recent article (\cite{MM20}), Mandal and Mishra have developed a theory of deformations of Hom-Lie-Rinehart algebras, which covers the case of Lie-Rinehart algebras. Our purpose here is to extend this theory to Lie-Rinehart superalgebras. We construct a suitable deformation complex and show that formal deformations are controlled by the cohomology obtained with this cochain complex. The case of Hom-Lie-Rinehart superalgebras will be studied in a forthcoming paper.

The article is organized as follows. In the first section, we recall the basics about superalgebras. Then, we give the definition of a Lie-Rinehart superalgebra in a second section, as well as some examples. The third section is devoted to the classification of Lie-Rinehart superalgebras, the associative superalgebra's dimension being less or equal than $2$ and the Lie superalgebra's dimension less or equal than $4$. We recall the classification of associative and Lie superalgebras and then we give all Lie-Rinehart superalgebra structures on all possible pairs. In Section 4, we construct a deformation cohomology complex and study formal deformations of Lie-Rinehart superalgebras. We show that the usual results about formal deformations remain true in this context. For example, we show that the infinitesimal element of such a deformation is a 2-cocycle of the complex constructed above. We also give a concrete example of deformation, relying on the classification that we made.

\section{Preliminaries}

Throughout this paper, we denote the group $\Z/n\Z$ by $\Z_n$, for $n\in \N$. If $V$ is a graded space, the degree of a homogeneous element $x\in V$ is denoted by $|x|$.

Let $\K$ be a commutative ring.
We summarize in the following the definitions of associative and Lie superalgebras, as well as related notions such as superderivations, using mainly \cite{KV77} or \cite{SM79}.

\begin{definition}
Let $V$ be a $\K$-module. It is said to be \textbf{$\Z_2$-graded} if it has a decomposition $V=V_0\oplus V_1$.
The elements of $V_i, i\in\Z_2$, are called homogeneous. An element $v\in V$ is called \textbf{even} if it belongs to $V_0$ and \textbf{odd} if it belongs to $V_1$. We denote $|v|=0$ if $v$ is even and $|v|=1$ if $v$ odd.
\end{definition}

\begin{definition}
An \textbf{associative $\K$-superalgebra} is a $\Z_2$-graded $\K$-module $A=A_0\oplus A_1$ endowed with a bilinear map $A\times A\longrightarrow A$ denoted by juxtaposition such that $(ab)c=a(bc)$ for all $a,b,c\in A$, and $A_iA_j\subset A_{i+j},$ the subscripts being taken modulo 2, \emph{i.e.} $i,j \in \Z_2$. 
\end{definition}

The gradings follow the rule $|ab|=|a|+|b|$ for all $a,b$ homogenous elements. The associativity of the bilinear map is the usual one, but the commutativity involves the gradings and is given by $ba=(-1)^{|a||b|}ab,$ for all homogenous elements $a,b$ in $A$. This identity is called \textbf{supercommutativity}.
We can sum this up by stating that the even elements commute with every other elements, and the odd ones anticommute with other odd elements. We then extend the supercommutativity to non-homogeneous elements. 

\begin{definition}\label{supermap}
Let $A$ and $B$ be two superalgebras. A \textbf{superalgebra morphism} is a linear map $f:A\longrightarrow B$ which satisfies the condition $f(ab)=(-1)^{|f||a|}f(a)f(b)$, for all $a,b\in A$. We say that $f$ is even and that the degree $|f|$ of $f$ is $0$ if $|f(a)|=|a|\mod(2)$ for all homogeneous elements $a\in A$, and we say that $f$ is odd and that and that the degree $|f|$ of $f$ is $1$ if $|f(a)|=(|a|+1)\mod(2)$. Thus, we always have the relation $|f(a)|=\left(|f|+|a| \right)\mod(2)$. 
\end{definition}

\begin{definition}\label{superderi}
Let $A$ be a superalgebra. A map $D:A\longrightarrow A$ is a \textbf{superderivation} (of degree
$|D|$) of $A$ if $D$ is a $\mathbb{Z}_2$-graded linear map and if the super-Leibniz condition is satisfied:
\[D(ab)=D(a)b+(-1)^{|a||D|}aD(b) \hspace{0.5cm} \forall a,b \in A. \]	
We denote by $\Der(A)$ the vector superspace of superderivations of $A$.
\end{definition}

\begin{definition}
A \textbf{Lie superalgebra} $L$ is a $\Z_2$-graded $\K$-module $L=L_0\bigoplus L_1$ endowed with a bracket $[\cdot,\cdot]$ satisfying, for all homogeneous elements $x,y,z\in L$:
\begin{enumerate}
\item $|[x,y]|=|x|+|y|$;
\item $[x,y]=-(-1)^{|x||y|}[y,x]$ (super-skewsymmetry);
\item $(-1)^{|x||z|}[x,[y,z]]+(-1)^{|z||y|}[z,[x,y]]+(-1)^{|x||y|}[y,[z,x]]=0$ (super-Jacobi).
\end{enumerate}	
\end{definition}

\begin{remark}
The super-Jacobi identity is equivalent to
\[
[x,[y,z]]=[[x,y],z]+(-1)^{|x||y|}[y,[x,z]]
\]
for all homogeneous element $x,y,z\in L$.
\end{remark}

\begin{example}
$\Der(A)$ carries a Lie superalgebra structure, with the bracket
\[ [D_1,D_2]=D_1\circ D_2-(-1)^{|D_1||D_2|}D_2\circ D_1.    \]
\end{example}

\begin{definition}
Let $L_1$ and $L_2$ be two Lie superalgebras. A \textbf{Lie superalgebra morphism} is a linear map $f:L_1 \to L_2$ which satisfies the condition
\[
f\left( [x,y]_{L_1}\right) =(-1)^{|f||x|}[f(x),f(y)]_{L_2},
\qquad \textup{for all } x,y\in L.
\]
We say that $f$ is even or odd following the same rule as in Definition \ref{supermap}.
\end{definition}

\begin{definition}
Let $V$ be a $\Z_2$-graded vector space and $A$ be a supercommutative superalgebra. We say that $V$ is a \textbf{left $\Z_2$-graded $A$-module} if there exists a map $A\times V\to V$, denoted by $(a,v)\mapsto av,$ such that $|av|=|a|+|v|$ and $a(bv)=(ab)v$, for all $v\in V$, $a,b\in A$.	
\end{definition}

\section{Lie-Rinehart superalgebras}

In this section, we focus on Lie-Rinehart superalgebras. We recall the basic definitions and give some examples, following Roger (\cite{RC20}) and Chemla (\cite{CS95}). 

\begin{definition}\label{superdef}
Let $\K$ be an arbitrary field. A \textbf{Lie-Rinehart superalgebra} is a pair $(A,L)$, where
\begin{itemize}
\item [$\bullet$] $L$ is a Lie superalgebra over $\K$, endowed with a bracket $[\cdot,\cdot]$;
\item [$\bullet$] $A$ is an associative and supercommutative $\K$-superalgebra,
\end{itemize}
such that, for all $x,y\in L$ and $a,b\in A$:
\begin{enumerate}
\item there is an action $A\times L\longrightarrow L,~(a,x)\longmapsto a\cdot x,$ making $L$ an $A$-module;

\item there is an action of $L$ on $A$ by superderivations: $L\to \text{Der}(A),~
x\to \left( \rho_x:a\mapsto \rho_x(a)\right)$, such that $\rho$ is an even morphism of Lie superalgebras;

\item $[x,ay]=\rho_x(a)y+(-1)^{|a||x|}a[x,y]$ (compatibility condition);

\item $\rho_{ax}(b)=a\rho_x(b)$ ($A$-linearity of $\rho$).
\end{enumerate}
\end{definition}

The maps in 1. and 2. must respect the gradings, \emph{i.e.}
\[
|ax|=|a|+|x|
\quad \textup{ and } \quad
|\rho_x(a)|=|a|+|x|.
\]
We sometimes write $\rho(x)(a)$ instead of $\rho_x(a)$.

\begin{remark}
If $A=\K$, the one dimensional even superalgebra, then any Lie-Rinehart superalgebra reduces in a certain sense to the ordinary Lie superalgebra over $\K$.
\end{remark}

\begin{example} 
Let $A$ be a supercommutative unital associative superalgebra and $L$ be a Lie superalgebra. Then the pair $(A,L)$ can always be endowed with a Lie-Rinehart superalgebra structure with the trivial action (the neutral element $e_0$ for the multiplication of $A$ acts by $e_0\cdot x=x$ for $x\in L$, and all the other elements of $A$ act by $0$) and the zero anchor ($\rho(x)=0~\forall x\in L$) (see Proposition \ref{trivialnull}).
\end{example}	

\begin{example}[Lie superalgebra of superderivations] Let $A$ be a supercommutative unital associative superalgebra, and $L=\Der(A)$ be its superalgebra of superderivations. Then one can check that the pair $\left( A, \Der(A)\right)$ is a Lie-Rinehart superalgebra, with the action $A\curvearrowright \Der(A)$ being given by $(a\cdot\delta)(b)=a\delta(b)$ and the trivial anchor being given by $\rho(\delta)=\delta$, for all $\delta \in \Der(A)$, $a,b\in A$.
\end{example}

\begin{example}[Crossed product] \label{supercrossed} We give now another example of Lie-Rinehart superalgebra, following Chemla (\cite{CS95}). Let $\mathfrak{g}$ be a Lie superalgebra, with bracket $[\cdot,\cdot]$, and $A$ a supercommutative associative superalgebra. We suppose that $\mathfrak{g}$ is endowed with a Lie superalgebra morphism $\sigma:\mathfrak{g}\longrightarrow\Der(A),~x\longmapsto\left( \sigma_x:a\longmapsto\sigma_x(a)\right).$

Then, $\sigma_x$ is a superderivation of $A$. We can define a new Lie superalgebra by setting $L:=A\otimes\mathfrak{g}$ endowed with the bracket
\[
[a\otimes x, b\otimes y]
=(-1)^{|x||b|}ab\otimes[x,y]+a\sigma_x(b)\otimes y-(-1)^{(|a|+|x|)(|b|+|y|)}b\sigma_y(a)\otimes x,
\]
with homogeneous $a,b\in A$ and homogeneous $x,y\in \mathfrak{g}$.

One can check that this formula gives a Lie superalgebra structure on $L$. We then extend $\sigma$ to a $A$-module morphism on $L=A\otimes\mathfrak{g}$, denoted by $\tilde{\sigma}\in\Der(A)$, in the following way: define	$ \tilde{\sigma}(a\otimes x):=a\sigma(x)=a\sigma_x$. We define a Lie-Rinehart structure on $(A,L)$, with $a,b,c\in A$ and $x,y\in \mathfrak{g}$:

\begin{enumerate}
\item Action $A\curvearrowright L$: $A\times L\longrightarrow L,~
(a,b\otimes y)\longmapsto a\cdot (b\otimes y)=ab\otimes y$,

\item Action by (super) derivations: $L\to \Der(A)
a\otimes x\mapsto \rho_{a\otimes x}:b\mapsto a\sigma_x(b) (=\tilde{\sigma}(a\otimes x)(b)).$
\end{enumerate}
\end{example}

\section{Classification in low dimensions}

In this section, we fix $\K=\C$, the complex field. We describe all Lie-Rinehart superstructures on a pair $(A,L)$ when $\dim(A)\leq 2$ and $\dim(L)\leq 4$. This means, if we are given a supercommutative associative superalgebra (with unit) $A$ and a Lie superalgebra $L$, we give all the pairs action-anchor which are compatible in the sense of Definition \ref{superdef}. 
We will denote basis elements of $A$ by $e_i^s$, $s=|e_i^s|\in\left\lbrace 0,1 \right\rbrace $, $1\leq i\leq \dim A_s$ and basis elements of $L$ by $(f_j^t)$, $t=| f_j^t|\in\left\lbrace 0,1 \right\rbrace $, $1\leq j\leq \dim L_t$.
Then we have
\[A_0=\left\langle e_i^0 \right\rangle_{{1\leq i\leq \dim A_0}},~ 
\A_1=\left\langle e_i^1 \right\rangle_{{1\leq i\leq \dim A_1}},~  
L_0=\left\langle f_j^0 \right\rangle_{{1\leq j\leq \dim L_0}},~  
L_1=\left\langle f_j^1 \right\rangle_{{1\leq j\leq \dim L_1}}.  \]

\subsection{Classification of associative and Lie superalgebras}

We recall the classification of supercommutative associative (\cite{AM14}, \cite{ACZ09} and \cite{GP74} for example) and Lie superalgebras (\cite{BN78} and \cite{BS99}) in low dimensions.

\subsubsection{Supercommutative associative superalgebras}\label{asso}

We list the supercommutative associative superalgebras with unit, of dimension up to 2. A list of all supercommutative associative superalgebras of dimension up to 4 can be found in Appendix \ref{appendixasso}. As above, we denote basis elements of $A_0$ by $e_i^0$ and those of $A_1$ by $e_j^1$. The unit is $e_1^0$. 
We only write the non zero products that have to be completed by supercommutativity and multilinearity. 

\begin{itemize}
\item[$\bullet$] The purely odd superalgebras $\bf{A_{0|p}}$ always have a zero product.
\item[$\bullet$] $\dim A=(1|0)$: there is only one superalgebra $\bf{A_{1|0}^1}$, whose product is given by $e_1^0e_1^0=e_1^0.$

\item[$\bullet$] $\dim A=(1|1)$: there is only one superalgebra $\bf{A_{1|1}^1}$ with product $e_1^1e_1^1=0$.

\item[$\bullet$] $\dim A=(2|0)$: there are two pairwise non-isomorphic superalgebras:
\begin{itemize}
\item [$\bf{A_{2|0}^1}:$] every non-trivial product is zero;
\item [$\bf{A_{2|0}^2}:$] $e_2^0e_2^0=e_2^0$.
\end{itemize}
\end{itemize}

\subsubsection{Lie superalgebras}\label{lialg}

We provide in the sequel a list of Lie superalgebras. As above, we denote the basis elements of $L_0$ by $f_{i}^0$ and those of $L_1$ by $f_{j}^1$. We only write the non zero brackets, the other brackets are obtained by super-skewsymmetry and bilinearity.

\begin{itemize}
\item[$\bullet$] The purely odd Lie superalgebras $\bf{L_{0|q}}$ always have a zero bracket.

\item[$\bullet$] $\dim L=(1|0)$: there is only one Lie superalgebra $\bf{L_{1|0}^1}$, whose bracket is given by $[f_1^0,f_1^0]=0.$

\item[$\bullet$] $\dim L=(1|1)$: there are three pairwise non-isomorphic Lie superalgebras:	
\begin{itemize}
\item [$\bf{L_{1|1}^1}:$] $[f_1^1,f_1^1]=f_1^0;$
\item [$\bf{L_{1|1}^2}:$] $[f_1^0,f_1^1]=f_1^1.$
\item [$\bf{L_{1|1}^3}:$] all brackets are null.
\end{itemize}

\item[$\bullet$] $\dim L=(1|2)$: there are six pairwise non-isomorphic Lie superalgebras:

\begin{itemize}
\item [$\bf{L_{1|2}^1}:$] $[f_1^0,f_1^1]=f_1^1,~ [f_1^0,f_2^1]=p f_2^1~(0<|p|\leq 1);$
\item [$\bf{L_{1|2}^2}:$] $[f_1^0,f_2^1]=f_1^1;$
\item [$\bf{L_{1|2}^3}:$] $[f_1^0,f_1^1]=f_1^1,~ [f_1^0,f_2^1]=f_1^1+f_2^1$;
\item [$\bf{L_{1|2}^4}:$] $[f_1^0,f_1^1]=p f_1^1-f_2^1,~ [f_1^0,f_2^1]=f_1^1+p f_2^1~(p\in\C);$
\item [$\bf{L_{1|2}^5}:$] $[f_1^1,f_1^1]=f_1^0,~ [f_2^1,f_2^1]=f_1^0$;
\item [$\bf{L_{1|2}^6}:$] null bracket.
\end{itemize}

\item[$\bullet$] $\dim L=(1|3)$: there are eight pairwise non-isomorphic Lie superalgebras:
\begin{itemize}
\item [$\bf{L_{1|3}^1}:$] $[f_1^0,f_1^1]=f_1^1,~ [f_1^0,f_2^1]=p f_2^1,~ [f_1^0,f_3^1]=q f_3^1~(0<|p|\leq|q|\leq 1);$
\item [$\bf{L_{1|3}^2}:$] $[f_1^0,f_1^1]=f_1^1,~[f_1^0,f_3^1]=f_2^1;$
\item [$\bf{L_{1|3}^3}:$] $[f_1^0,f_1^1]=pf_1^1,~ [f_1^0,f_2^1]=f_2^1,~ [f_1^0,f_3^1]=f_2^1+f_3^1$;
\item [$\bf{L_{1|3}^4}:$] $[f_1^0,f_1^1]=pf_1^1,~ [f_1^0,f_2^1]=qf_2^1-f_3^1,~ [f_1^0,f_3^1]=f_2^1+qf_3^1,~p\neq 0$;
\item [$\bf{L_{1|3}^5}:$] $[f_1^0,f_2^1]=f_1^1,~ [f_1^0,f_3^1]=f_2^1,~p,q\neq 0$;
\item [$\bf{L_{1|3}^6}:$] $[f_1^0,f_1^1]=f_1^1,~ [f_1^0,f_2^1]=f_1^1+f_2^1,~ [f_1^0,f_3^1]=f_2^1+f_3^1$;
\item [$\bf{L_{1|3}^7}:$] $[f_1^1,f_1^1]=f_1^0,~ [f_2^1,f_2^1]=f_1^0,~ [f_3^1,f_3^1]=f_1^0$;
\item [$\bf{L_{1|3}^8}:$] null bracket.
\end{itemize}

\item[$\bullet$] $\dim L=(2|0)$: there are two pairwise non-isomorphic Lie superalgebras:
\begin{itemize}
\item [$\bf{L_{2|0}^1}:$] null bracket;
\item [$\bf{L_{2|0}^2}:$] $[f_1^0,f_2^0]=f_2^0.$
\end{itemize}

\item[$\bullet$] $\dim L=(2|1)$: there are six pairwise non-isomorphic Lie superalgebras:
\begin{itemize}
\item [$\bf{L_{2|1}^1}:$] $[f_1^1,f_1^1]=f_2^0;$
\item [$\bf{L_{2|1}^2}:$] $[f_1^0,f_1^1]=f_1^1~[f_2^0,f_1^1]=-f_1^1;$
\item [$\bf{L_{2|1}^3}:$] $[f_1^0,f_2^0]=f_2^0,~ [f_1^0,f_1^1]=\frac{1}{2}f_1^1$;
\item [$\bf{L_{2|1}^4}:$] $[f_1^0,f_2^0]=f_2^0,~ [f_1^0,f_1^1]=p f_1^1,~(p \neq 0)$;
\item [$\bf{L_{2|1}^5}:$] $[f_1^0,f_2^0]=f_2^0$.
\item [$\bf{L_{2|1}^6}:$] null bracket.
\end{itemize}

\item[$\bullet$] $\dim L=(2|2)$: there are eighteen pairwise non-isomorphic Lie superalgebras:
\begin{itemize}
\item [$\bf{L_{2|2}^1}:$] $[f_1^0,f_1^1]=f_1^1,~[f_1^0,f_2^1]=f_2^1,~[f_2^0,f_2^1]=f_1^1;$
\item [$\bf{L_{2|2}^2}:$] $[f_1^0,f_1^1]=f_1^1,~[f_1^0,f_2^1]=f_2^1,~[f_2^0,f_2^1]=f_1^1,~[f_2^0,f_1^1]=-f_2^1;$
\item [$\bf{L_{2|2}^3}:$]$[f_1^0,f_2^0]=f_2^0,~[f_1^0,f_1^1]=pf_1^1,~[f_1^0,f_2^1]=qf_2^1,~pq\neq0;$
\item [$\bf{L_{2|2}^4}:$]$[f_1^0,f_2^0]=f_2^0,~[f_1^0,f_1^1]=pf_1^1,~[f_1^0,f_2^1]=f_1^1+pf_2^1,~p\neq 0;$ 
\item [$\bf{L_{2|2}^5}:$] $[f_1^0,f_2^0]=f_2^0,~[f_1^0,f_1^1]=pf_1^1-qf_2^1,~[f_1^0,f_2^1]=qf_1^1+pf_2^1,~q\neq 0;$
\item [$\bf{L_{2|2}^6}:$] $[f_1^0,f_2^0]=f_2^0,~[f_1^0,f_1^1]=(p+1)f_1^1,~[f_1^0,f_2^1]=pf_2^1,~[f_2^0,f_2^1]=f_1^1,~p\neq 0;$
\item [$\bf{L_{2|2}^7}:$] $[f_1^0,f_2^0]=f_2^0,~[f_1^0,f_1^1]=\frac{1}{2}f_1^1,~[f_1^0,f_2^1]=\frac{1}{2}f_2^1,~[f_1^1,f_1^1]=f_2^0,~[f_2^1,f_2^1]=f_2^0;$
\item [$\bf{L_{2|2}^8}:$] $[f_1^0,f_2^0]=f_2^0,~[f_1^0,f_1^1]=\frac{1}{2}f_1^1,~[f_1^0,f_2^1]=\frac{1}{2}f_2^1,~[f_1^1,f_1^1]=f_2^0;$
\item [$\bf{L_{2|2}^9}:$] $[f_1^0,f_2^0]=f_2^0,~[f_1^0,f_1^1]=pf_1^1,~[f_1^0,f_2^1]=(1-p)f_2^1,~[f_1^1,f_2^1]=f_2^0;$
\item [$\bf{L_{2|2}^{10}}:$] $[f_1^0,f_2^0]=f_2^0,~[f_1^0,f_1^1]=\frac{1}{2}f_1^1,~[f_1^0,f_2^1]=f_1^1+\frac{1}{2}f_2^1,~[f_2^1,f_2^1]=f_2^0;$
\item [$\bf{L_{2|2}^{11}}:$] $[f_1^0,f_2^0]=f_2^0,~[f_1^0,f_1^1]=\frac{1}{2}f_1^1-pf_2^1,~[f_1^0,f_2^1]=pf_1^1+\frac{1}{2}f_2^1,~[f_1^1,f_1^1]=f_2^0$ and $[f_2^1,f_2^1]=f_2^0$, for $p\neq0$;
\item [$\bf{L_{2|2}^{12}}:$] $[f_1^0,f_2^0]=f_2^0,~[f_1^0,f_1^1]=f_1^1,~[f_2^0,f_2^1]=f_1^1,~[f_1^1,f_2^1]=-\frac{1}{2}f_2^0,~[f_2^1,f_2^1]=f_1^0;$
\item [$\bf{L_{2|2}^{13}}:$] $[f_1^1,f_1^1]=f_1^0,~[f_2^1,f_2^1]=f_2^0;$
\item [$\bf{L_{2|2}^{14}}:$] $[f_1^0,f_2^0]=f_2^0,~[f_1^0,f_1^1]=f_1^1,~[f_1^1,f_2^1]=f_2^0;$
\item [$\bf{L_{2|2}^{15}}:$] $[f_1^0,f_2^0]=f_2^0,~[f_1^0,f_1^1]=\frac{1}{2}f_1^1,~[f_1^1,f_1^1]=f_2^0;$
\item [$\bf{L_{2|2}^{16}}:$] $[f_1^0,f_1^1]=f_1^1,~[f_1^0,f_{2}^1]=-f_2^1,~[f_1^1,f_2^1]=f_2^0;$
\item [$\bf{L_{2|2}^{17}}:$] $[f_1^0,f_2^1]=f_1^1,~[f_2^1,f_2^1]=f_2^0;$
\item [$\bf{L_{2|2}^{18}}:$] null bracket.
\end{itemize}

\item[$\bullet$] $\dim L=(3|0)$: there are six pairwise non-isomorphic Lie superalgebras:
\begin{itemize}
\item [$\bf{L_{3|0}^1}:$] null bracket;
\item [$\bf{L_{3|0}^2}:$] $[f_1^0,f_2^0]=f_3^0$;
\item [$\bf{L_{3|0}^3}:$] $[f_1^0,f_2^0]=f_1^0$;
\item [$\bf{L_{3|0}^4}:$] $[f_1^0,f_2^0]=f_2^0$, $[f_1^0,f_3^0]=f_2^0+f_3^0$;
\item [$\bf{L_{3|0}^5}:$] $[f_1^0,f_2^0]=f_2^0$, $[f_1^0,f_3^0]=pf_3^0$, $p \neq 0$;
\item [$\bf{L_{3|0}^6}:$] $[f_1^0,f_2^0]=f_3^0$, $[f_1^0,f_3^0]=-2f_1^0$, $[f_2^0,f_3^0]=2f_2^0$;
\end{itemize}

\item[$\bullet$] $\dim L=(3|1)$: there are seven pairwise non-isomorphic Lie superalgebras:
\begin{itemize}
\item [$\bf{L_{3|1}^1}:$] $[f_2^0,f_3^0]=f_1^0,~[f_2^0,f_1^1]=f_1^1$
\item [$\bf{L_{3|1}^2}:$] $[f_1^0,f_3^0]=f_1^0,~[f_2^0,f_3^0]=f_1^0+f_2^0,~[f_3^0,f_1^1]=qf_1^1,~q\neq 0;$
\item [$\bf{L_{3|1}^3}:$] $[f_1^0,f_3^0]=pf_1^0-f_2^0,~[f_2^0,f_3^0]=f_1^0+pf_2^0,~[f_3^0,f_1^1]=qf_1^1,~pq\neq 0;$
\item [$\bf{L_{3|1}^4}:$] $[f_2^0,f_3^0]=f_1^0,~ [f_1^1,f_1^1]=f_1^0$;
\item [$\bf{L_{3|1}^5}:$] $[f_1^0,f_2^0]=f_2^0,~[f_1^0,f_3^0]=pf_3^0,~[f_1^0,f_1^1]=\frac{1}{2}f_1^1,~[f_1^1,f_1^1]=f_2^0,~p\neq 0$;
\item [$\bf{L_{3|1}^6}:$] $[f_1^0,f_2^0]=f_2^0,~[f_1^0,f_3^0]=f_2^0+f_3^0,~[f_1^0,f_1^1]=\frac{1}{2}f_1^1,~[f_1^1,f_1^1]=f_2^0$.
\item [$\bf{L_{3|1}^7}:$] null bracket.
\end{itemize}

\item[$\bullet$] $\dim L=(4|0)$: there are sixteen pairwise non-isomorphic Lie superalgebras:
\begin{itemize}
\item [$\bf{L_{4|0}^1}:$] null bracket;
\item [$\bf{L_{4|0}^2}:$] $[f_1^0,f_2^0]=f_3^0$;
\item [$\bf{L_{4|0}^3}:$] $[f_1^0,f_2^0]=f_1^0$;
\item [$\bf{L_{4|0}^4}:$] $[f_1^0,f_2^0]=f_2^0$, $[f_1^0,f_3^0]=f_2^0+f_3^0$;
\item [$\bf{L_{4|0}^5}:$] $[f_1^0,f_2^0]=f_2^0$, $[f_1^0,f_3^0]=pf_3^0$, $0<|p|\leq 1 $;
\item [$\bf{L_{4|0}^6}:$] $[f_1^0,f_2^0]=f_1^0$, $[f_3^0,f_4^0]=f_3^0$;
\item [$\bf{L_{4|0}^7}:$] $[f_1^0,f_2^0]=f_3^0$, $[f_1^0,f_3^0]=-2f_1^0$, $[f_2^0,f_3^0]=2f_2^0$;
\item [$\bf{L_{4|0}^8}:$] $[f_1^0,f_2^0]=f_3^0$, $[f_1^0,f_3^0]=f_4^0$;
\item [$\bf{L_{4|0}^9}:$] $[f_1^0,f_2^0]=f_2^0$, $[f_1^0,f_3^0]=f_3^0$, $[f_1^0,f_4^0]=pf_3^0$, $p\neq 0$;
\item [$\bf{L_{4|0}^{10}}:$] $[f_1^0,f_2^0]=f_3^0$, $[f_1^0,f_3^0]=f_4^0$, $[f_1^0,f_4^0]=pf_2^0-qf_3^0+f_4^0$, $(p,q)\in\C^{\times}\times\C$ or $(p,q)=(0,0)$;
\item [$\bf{L_{4|0}^{11}}:$] $[f_1^0,f_2^0]=f_3^0$, $[f_1^0,f_3^0]=f_4^0$, $[f_1^0,f_4^0]=p(f_2^0+f_3^0)$, $p\neq 0$;
\item [$\bf{L_{4|0}^{12}}:$] $[f_1^0,f_2^0]=f_3^0$, $[f_1^0,f_3^0]=f_4^0$, $[f_1^0,f_4^0]=f_2^0$;
\item [$\bf{L_{4|0}^{13}}:$] $[f_1^0,f_2^0]=\frac{1}{3}f_2^0+f_3^0$, $[f_1^0,f_3^0]=\frac{1}{3}f_3^0$, $[f_1^0,f_4^0]=\frac{1}{3}f_4^0$
\item [$\bf{L_{4|0}^{14}}:$] $[f_1^0,f_2^0]=f_2^0$, $[f_1^0,f_3^0]=f_3^0$, $[f_1^0,f_4^0]=2f_4^0$, $[f_2^0,f_3^0]=f_4^0$;
\item [$\bf{L_{4|0}^{15}}:$] $[f_1^0,f_2^0]=f_3^0$, $[f_1^0,f_3^0]=f_2^0$, $[f_2^0,f_3^0]=f_4^0$;
\item [$\bf{L_{4|0}^{16}}:$] $[f_1^0,f_2^0]=f_3^0$, $[f_1^0,f_3^0]=-pf_2^0+f_3^0$, $[f_1^0,f_4^0]=f_4^0$, $[f_2^0,f_3^0]=f_4^0,$ $p\neq 0$;
\end{itemize}
\end{itemize}

\subsection{Superderivations of associative superalgebras}

If $A$ is a superalgebra, we have given the definition of superderivations of $A$ in Definition~\ref{superderi}. Let $(A,L)$ be a Lie-Rinehart superalgebra. As an anchor $\rho: L\longrightarrow L$ corresponds to a family $\left\lbrace \rho(x) \right\rbrace_{x\in L}$ of superderivations of $A$, it seems quite natural to study and describe  superderivations spaces for all the superalgebras we deal with. 
We recall that in this case, $\Der(A)$ has a Lie superalgebra structure, the bracket being given by the supercommutator.

For every associative superalgebra listed above, we give the general form of the superderivations $D$. We describe $D$ with respect to the basis $\left\{e_1^0, e_2^0, \dotsc, e_n^0, e_1^1, e_2^1, \dotsc, e_p^1 \right\}$ of the algebra $\bf{A_{n|p}^k}$. We only write the non-zero values. All the parameters are complex and independent in each column.

\noindent\noindent
\resizebox{\textwidth}{!}{%
\parbox{1.23\textwidth}{%
\begin{tabular}{|c|c|c|}
\hline
\textbf{Superalgebra} & \textbf{Superderivations of degree $0$} & \textbf{Superderivations of degree $1$} \\
\hline
$\bf{A_{1|0}^1}$&$0$&$0$\\
\hline
$\bf{A_{1|1}^1}$&$D(e_1^1)=\lambda_1e_1^1$&$D(e_1^1)=\lambda_1e_1^0$\\
\hline
$\bf{A_{1|2}^1}$&$D(e_1^1)=\lambda_1e_1^1+\lambda_2e_2^1,~D(e_2^1)=\lambda_3e_1^1+\lambda_4e_2^1 $&null\\
\hline
&$D(e_1^1)=\lambda_1e_1^1+\lambda_2e_2^1+\lambda_3e_3^1$&\\$\bf{A_{1|3}^1}$&$D(e_2^1)=\lambda_4e_1^1+\lambda_5e_2^1+\lambda_6e_3^1$&null\\&$D(e_3^1)=\lambda_7e_1^1+\lambda_8e_2^1+\lambda_9e_3^1$&\\
\hline
$\bf{A_{2|0}^1}$&$D(e_2^0)=\lambda_1e_2^0$&null\\
\hline
$\bf{A_{2|0}^2}$&null&null\\
\hline
$\bf{A_{2|1}^1}$&$D(e_2^0)=\lambda_1e_2^0,~D(e_1^1)=\lambda_2e_1^1$&$D(e_2^0)=\lambda_1e_1^1,~D(e_1^1)=\lambda_2e_2^0$\\
\hline
$\bf{A_{2|1}^2}$&$D(e_1^1)=\lambda_1e_1^1$&$D(e_1^1)=\lambda_1e_2^0$\\
\hline
$\bf{A_{2|2}^1}$&$D(e_1^1)=\lambda_1e_1^1,~D(e_2^1)=\lambda_2e_2^1$&$D(e_1^1)=\lambda_1e_2^0,~D(e_2^1)=\lambda_2(e_1^0-e_2^0)$\\
\hline
$\bf{A_{2|2}^2}$&$D(e_1^1)=\lambda_1e_1^1+\lambda_2e_2^1,~D(e_2^1)=\lambda_3e_1^1+\lambda_4e_2^1$&null\\
\hline
\multirow{2}{*}{$\bf{A_{2|2}^3}$}&$D(e_2^0)=\lambda_1e_2^0,~D(e_1^1)=\lambda_2e_1^1+\lambda_3e_2^1,$&$D(e_2^0)=\lambda_1e_2^1,~D(e_1^1)=\lambda_2e_1^0+\lambda_3e_2^0,$\\&$D(e_2^1)=(\lambda_1+\lambda_2)e_2^1$&$D(e_2^1)=\lambda_2e_2^0$\\
\hline
\multirow{2}{*}{$\bf{A_{2|2}^4}$}&$D(e_2^0)=\lambda_1e_2^0,~D(e_1^1)=\lambda_2e_1^1+\lambda_3e_2^1,$&$D(e_2^0)=\lambda_1e_1^1+\lambda_2e_2^1,~D(e_1^1)=\lambda_3e_2^0$\\&$D(e_2^1)=\lambda_4e_1^1+\lambda_5e_2^1$&$D(e_2^1)=\lambda_2e_2^0$\\
\hline
\multirow{2}{*}{$\bf{A_{2|2}^5}$}&$D(e_2^0)=\lambda_1e_2^0,~D(e_1^1)=\lambda_2e_1^1+\lambda_3e_2^1,$&$D(e_2^0)=\lambda_1e_1^1+\lambda_2e_2^1,~D(e_1^1)=\lambda_2e_1^0+\lambda_3e_2^0$\\&$D(e_2^1)=\lambda_4e_1^1+(\lambda_1-\lambda_2)e_2^1$&$D(e_2^1)=-\lambda_1e_1^0+\lambda_4e_2^0$\\
\hline
$\bf{A_{3|0}^1}$&null&null\\
\hline
$\bf{A_{3|0}^2}$&$D(e_3^0)=\lambda_1e_3^0$&null\\
\hline
$\bf{A_{3|0}^3}$&$D(e_3^0)=\lambda_1e_3^0$&null\\
\hline
$\bf{A_{3|0}^4}$&$D(e_2^0)=\lambda_1e_2^0+\lambda_2e_3^0,~D(e_3^0)=\lambda_3e_2^0+\lambda_4e_3^0$&null\\
\hline
$\bf{A_{3|1}^1}$&$D(e_1^1)=\lambda_1e_1^1$&$D(e_1^1)=\lambda_1(e_1^0-e_2^0-e_3^0)$\\
\hline
$\bf{A_{3|1}^2}$&$D(e_3^0)=\lambda_1e_3^0,~D(e_1^1)=\lambda_2e_1^1$&$D(e_1^1)=\lambda_1(e_1^0-e_2^0)$\\
\hline
$\bf{A_{3|1}^3}$&$D(e_3^0)=\lambda_1e_3^0,~D(e_1^1)=\lambda_2e_1^1$&$D(e_3^0)=\lambda_1e_1^1,~D(e_1^1)=\lambda_2e_3^0$\\
\hline
\multirow{2}{*}{$\bf{A_{3|1}^4}$}&$D(e_2^0)=\lambda_1e_2^0+\lambda_2e_3^0,~ D(e_3^0)=2\lambda_1e_3^0, $&	\multirow{2}{*}{$D(e_2^0)=\lambda_1e_1^1,~D(e_1^1)=\lambda_2e_3^0$}\\&$D(e_1^1)=\lambda_3e_1^1$&\\
\hline
\multirow{2}{*}{$\bf{A_{3|1}^5}$}&$D(e_2^0)=\lambda_1e_2^0+\lambda_2e_3^0,~D(e_3^0)=\lambda_3e_2^0+\lambda_4e_3^0 $&$D(e_2^0)=\lambda_1e_1^1,~D(e_3^0)=\lambda_2e_1^1,$ \\&$D(e_1^1)=\lambda_5e_1^1$&$D(e_1^1)=\lambda_3e_2^0+\lambda_4e_3^0$\\
\hline
$\bf{A_{4|0}^1}$&null&null\\
\hline
$\bf{A_{4|0}^2}$&$D(e_4^0)=\lambda_1e_4^0$&null\\
\hline
$\bf{A_{4|0}^3}$&$D(e_3^0)=\lambda_1e_3^0,~D(e_4^0)=\lambda_2e_4^0$&null\\
\hline
$\bf{A_{4|0}^4}$&$D(e_3^0)=\lambda_1e_3^0+\lambda_2e_4^0,~D(e_4^0)=2\lambda_1e_4^0$&null\\
\hline
\multirow{2}{*}{$\bf{A_{4|0}^5}$}&$D(e_2^0)=\lambda_1e_2^0+\lambda_2e_3^0+\lambda_3e_4^0,~D(e_4^0)=3\lambda_1e_4^0,$&\multirow{2}{*}{null}\\&$D(e_3^0)=2\lambda_1e_3^0+2\lambda_2e_4^0$&\\
\hline
$\bf{A_{4|0}^6}$&$D(e_3^0)=\lambda_1e_3^0+\lambda_2e_4^0,~D(e_4^0)=\lambda_3e_3^0+\lambda_4e_4^0$&null\\
\hline
$\bf{A_{4|0}^7}$&$D(e_2^0)=\lambda_1e_2^0+\lambda_2e_4^0,~D(e_4^0)=\lambda_3e_3^0+\lambda_24e_4^0$&null\\
\hline
\multirow{2}{*}{$\bf{A_{4|0}^8}$}&$D(e_2^0)=\lambda_1e_2^0+\lambda_2e_3^0+\lambda_3e_4^0,~D(e_3^0)=2\lambda_1e_3^0$&\multirow{2}{*}{null}\\&$D(e_4^0)=\lambda_4e_3^0+\lambda_5e_4^0$&\\
\hline
&$D(e_2^0)=\lambda_1e_2^0+\lambda_2e_3^0+\lambda_3e_4^0$&\\$\bf{A_{4|0}^9}$&$D(e_3^0)=\lambda_4e_2^0+\lambda_5e_3^0+\lambda_6e_4^0$&null\\&$D(e_4^0)=\lambda_7e_2^0+\lambda_8e_3^0+\lambda_9e_4^0$&\\
\hline
\end{tabular}
}}

\subsection{Classification of Lie-Rinehart superalgebras}

We provide in this section a classification of Lie-Rinehart superalgebras in low dimensions, using the following general results and the computer algebra system Mathematica. We write only non trivial and non zero relations. If $(A,L)$ is a Lie-Rinehart superalgebra, we denote its dimension by a tuple $(n|p,m|q)$, where $n=\dim A_0$, $p=\dim A_1$, $m=\dim L_0$, $q=\dim L_1$. We say that an action is \textbf{trivial} if $e_i^s\cdot f_j^t=f_j^t$ if $i=1, s=0$, and $0$ otherwise. Using properties of the degrees and basic calculations, we obtain some general results.

\begin{proposition}\label{trivialnull}
If $A$ is a supercommutative associative superalgebra and $L$ a Lie superalgebra, then we can always endow the pair $(A,L)$ with a Lie-Rinehart superstructure using the trivial action and the null anchor.
\end{proposition}

\begin{proof}
If $\rho(x)(a)=0~\forall x\in L,~\forall a\in A$, it's clear that $\rho:L\longrightarrow\Der(A)$ is an $A$-linear morphism of Lie superalgebras and that $\rho(x)$ is a superderivation. For $y\in L$, the compatibility condition is then
\[ [x,ay]=(-1)^{|a||x|}a[x,y], \] which is always satisfied, because the action is trivial. 
\end{proof}

\begin{proposition} Let $(A,L)$ be a Lie-Rinehart superalgebra. If $\dim(A,L)=(n|p,0|0)$, $(1|0,m|q)$, or $(1|p,0|q)$, then the only possible structure is given by the trivial action and the null anchor.
\end{proposition}
\begin{proof}
\begin{enumerate} 	
\item $(n|p,0|0)$ case: immediate;
\item $(1|0,m|q)$ case: if $e_1^0$ is the unit of $A$, then $e_1^0\cdot x=x$ $\forall x\in L$ and $\rho(x)(e_1^0)=0$ $\forall x\in L$;
\item $(1|p,0|q)$ case: $L_0=\left\lbrace0\right\rbrace$, so the action is trivial. Let $1\leq k\leq q$ and $1\leq l\leq p$. We have \[\rho(f_k^1)(e_l^1)=r_{(k,1)(l,1)}^1e_1^0,~r_{(k,1)(l,1)}^1\in \C .\] Then, using the fourth condition of the Definition \ref{superdef},
\[0=\rho(e_l^1\cdot f_k^1)(e_l^1) =e_l^1\rho(f_k^1)(e_l^1)=e_1^1\cdot r_{(k,1)(l,1)}^1e_1^0=r_{(k,1)(l,1)}^1e_l^1. \]
So $r_{(k,1)(l,1)}^1=0$ and the anchor vanishes.
\qedhere
\end{enumerate}
\end{proof}

There are also some exceptional cases, found by computer calculations, where the only suitable pair is the trivial action and the zero anchor, given by the following proposition:

\begin{proposition} For the following Lie-Rinehart superalgebras, the only compatible action/anchor pair is the trivial action and the zero anchor:
\[ 
\bf{\left( A_{1|1}^1, L_{2|2}^{12} \right)};~ 
\bf{\left( A_{1|1}^1, L_{3|0}^{6} \right)}; ~ 
\bf{\left( A_{2|0}^1, L_{1|1}^{1} \right)}; ~ 
\bf{\left( A_{2|0}^1, L_{1|2}^{2} \right)}; ~ 
\bf{\left( A_{2|0}^1, L_{1|3}^{7} \right)}.
\]
\end{proposition}

One notices that all supercommutative associative superalgebras in the above list have a trivial multiplication, $i.e.$ all the products vanishes, except the products involving the unit. One may conjecture that if the only suitable pair of a Lie-Rinehart superalgebra $(A,L)$ is the trivial action and the zero anchor, then $A$ must have a trivial multiplication. This result is to be proven yet.

\begin{proposition}
\begin{itemize}
\item [$\bullet$] If $ \dim(A,L)=(0|p,m|0)$, the action vanishes;
\item [$\bullet$] If $ \dim(A,L)=(n|0,0|q)$, the anchor vanishes.
\item [$\bullet$] If $ \dim(A,L)=(n|p,m|0)$ or $ (A,L)=(n|p,0|q)$, the elements of $A_1$ are acting by $0$;
\item [$\bullet$] If $ \dim(A,L)=(n|0,m|q)$, 
the elements of $L_1$ are acting (via the anchor) by $0$.		
\end{itemize}
\end{proposition}

Now, we list the Lie-Rinehart superalgebras which are not already classified by the above results. They are arranged in lexicographic order of the tuple $(n|p,m|q)$. For each tuple, a table gives all the possible pairs, with all compatible actions and anchors. Every row of a table gives a different Lie-Rinehart superstructure. We give here all the tables with $\dim(A)\leq 2$ and $\dim(L)\leq2$. The tables with $\dim(L)>2$ are given in Appendix \ref{appendix2}.

\begin{description}
\item[$(1|1, 1|0)$-type:]
The only possible Lie-Rinehart superalgebra is given by the trivial action and the anchor $\rho(f_1^0)(e_1^1)=\lambda e_1^1,$ $\lambda\in\C$.

\item[$(1|1, 1|1)$-type:] $(\lambda\in\C)$	
\begin{center}
\begin{tabular}{|c|c|l|l|}
\hline
$\textbf{A}$ & $\textbf{L}$ & $\textbf{Action}$ & $\textbf{Anchor}$ \\
\hline
\multirow{7}{*}{$A_{1|1}^1$} & $L_{1|1}^1$ & $e_1^1\cdot f_1^1=\lambda f_1^0$ & null\\
\cline{2-4}
& \multirow{3}{*}{$L_{1|1}^2$}& trivial & $\rho(f_1^0)(e_1^1)=\lambda e_1^1$ \\
\cline{3-4}
& & $e_1^1\cdot f_1^1=\lambda f_1^0$& $\rho(f_1^0)(e_1^1)=-e_1^1,~\rho(f_1^1)(e_1^1)=-\lambda e_1^0$ \\
\cline{3-4}
& & $e_1^1\cdot f_1^0=\lambda f_1^1$ &$\rho(f_1^0)(e_1^1)=e_1^1$ \\
\cline{2-4}
& \multirow{2}{*}{$L_{1|1}^3$}& trivial & $\rho(f_1^0)(e_1^1)=\lambda e_1^1$ \\
\cline{3-4}
& & $e_1^1\cdot f_1^0=\lambda f_1^1$ & null \\
\hline
\end{tabular}
\end{center}

\item[$(1|1, 2|0)$-type:] $(\lambda, \mu\in\C)$
\begin{center}
\begin{tabular}{|c|c|l|l|}
\hline
$\textbf{A}$ & $\textbf{L}$ & $\textbf{Action}$ & $\textbf{Anchor}$ \\
\hline
\multirow{2}{*}{$A_{1|1}^1$} & $L_{2|0}^1$ & trivial & $\rho(f_1^0)(e_1^1)=\lambda e_1^1,~\rho(f_2^0)(e_1^1)=\mu e_1^1$ \\
\cline{2-4}
& $L_{2|0}^2$ & trivial & $\rho(f_1^0)(e_1^1)=\lambda e_1^1$ \\
\hline
\end{tabular}
\end{center}

\item[$(2|0, 0|1)$-type:] The only remarkable pair is $(\mathbf{A_{2|0}^2}, \mathbf{L_{0|1}^1})$, endowed with the null anchor and the action $e_2^0\cdot f_1^1=f_1^1$.

\item[$(2|0, 0|2)$-type:]
Here $L=\bf{L_{0|2}^1}$ and the anchor is always null. We list the possible compatible actions for each supercommutative associative $(2|0)$-type superalgebra:
\begin{itemize}
\item[$\bullet$] $A=A_{2|0}^1$:
\begin{enumerate}
\item $e_2^0\cdot f_2^1=\lambda f_1^1,~\lambda \in\C$;
\item $e_2^0\cdot f_1^1=\lambda f_1^1+\mu f_2^1$, $e_2^0\cdot f_2^1=-\frac{\lambda^2}{\mu}f_1^1-\mu f_2^1$, $(\lambda,\mu)\in \C\times\C^{\times}$.
\end{enumerate}
\item[$\bullet$] $A=A_{2|0}^2$:
\begin{enumerate}
\item $e_2^0\cdot f_2^1=\lambda f_1^1+f_2^1,~\lambda \in\C$;
\item $e_2^0\cdot f_1^1=f_1^1$, $e_2^0\cdot f_2^1=\lambda f_1^1$, $\lambda\in\C$;
\item $e_2^0\cdot f_2^1=f_2^1$;
\item $e_2^0\cdot f_1^1=f_1^1$, $e_2^0\cdot f_2^1=f_2^1$;
\item $e_2^0\cdot f_1^1=\lambda f_1^1+\mu f_2^1$, $e_2^0\cdot f_2^1=\frac{\lambda-\lambda^2}{\mu}f_1^1-(1-\mu)f_2^1$, $(\lambda,\mu)\in \C\times\C^{\times}$.
\end{enumerate}	
\end{itemize}

\item[$(2|0, 1|0)$-type:] $(\lambda\in\C)$
\begin{center}
\begin{tabular}{|c|c|c|c|}
\hline
$\textbf{A}$ & $\textbf{L}$ & $\textbf{Action}$ & $\textbf{Anchor}$ \\
\hline
$A_{2|0}^1$ & $L_{1|0}$ & trivial & $\rho(f_1^0)(e_2^0)=\lambda e_2^0$ \\ 
\hline
$A_{2|0}^2$ & $L_{1|0}$ & $e_2^0\cdot f_1^0=f_1^0$ & null \\
\hline
\end{tabular}
\end{center}

\item[$(2|0, 1|1)$-type:] $(\lambda\in\C)$
\begin{center}
\begin{tabular}{|c|c|c|c|}
\hline
$\textbf{A}$ & $\textbf{L}$ & $\textbf{Action}$ & $\textbf{Anchor}$ \\
\hline
\multirow{2}{*}{$A_{2|0}^1$} & $L_{1|1}^2$ & trivial & $\rho(f_1^0)(e_2^0)=\lambda e_2^0$ \\
\cline{2-4}
& $L_{1|1}^3$ & trivial & $\rho(f_1^0)(e_2^0)=\lambda e_2^0$ \\
\hline
\multirow{5}{*}{$A_{2|0}^2$} & \multirow{1}{*}{\smash{$L_{1|1}^1~\&~L_{1|1}^2$}} & $e_2^0\cdot f_1^0=f_1^0,~e_2^0\cdot f_1^1=f_1^1$ & \multirow{1}{*}{null}\\
\cline{2-4}
& \multirow{3}{*}{$L_{1|1}^3$} & $e_2^0\cdot f_1^1=f_1^1$ & \multirow{3}{*}{null} \\
\cline{3-3}
& & $e_2^0\cdot f_1^0=f_1^0$ &\\
\cline{3-3}
& & $e_2^0\cdot f_1^0=f_1^0,~ e_2^0\cdot f_1^1=f_1^1 $ &\\
\hline
\end{tabular}
\end{center}

\begin{remark}
We see here that for $A=A_{2|0}^2$, all the anchors are null. It is a straightforward consequence of the fact that the super derivations space of $A_{2|0}^2$ is $\left\lbrace 0\right\rbrace$. 
\end{remark}

\item[$(2|0, 2|0)$-type:] $(\lambda\in\C)$

\noindent
\resizebox{\textwidth}{!}{%
\parbox{1.2\textwidth}{%
\begin{tabular}{|c|c|l|c|}
\hline
$\textbf{A}$ & $\textbf{L}$ & $\textbf{Action}$ & $\textbf{Anchor}$ \\[1pt]
\hline
\multirow{5}{*}{$A_{2|0}^1$} & \multirow{3}{*}{$L_{2|0}^1$} & trivial & $\rho(f_1^0)(e_2^0)=\lambda e_2^0,~\rho(f_2^0)(e_2^0)=\mu e_2^0$ \\[1pt]
\cline{3-4}
&& $e_2^0\cdot f_2^0=\lambda f_1^0$ & null \\[1pt]
\cline{3-4}
&& $e_2^0\cdot f_1^0=\lambda f_1^0+\mu f_2^0,~e_2^0\cdot f_2^0=-\frac{\lambda^2}{\mu}f_1^0-\lambda f_2^0,~\mu \neq0$ & null \\[1pt]
\cline{2-4}
& \multirow{2}{*}{$L_{2|0}^2$} & trivial & $\rho(f_1^0)(e_2^0)=\lambda e_2^0$ \\[1pt]
\cline{3-4}
&& $e_2^0\cdot f_1^0=\lambda f_2^0$ & $\rho(f_1^0)(e_2^0)=e_2^0$ \\[1pt]
\hline
\multirow{7}{*}{$A_{2|0}^2$} & \multirow{6}{*}{$L_{2|0}^1$} & $e_2^0\cdot f_1^0= f_1^0,~e_2^0\cdot f_2^0= f_2^0$ & \multirow{6}{*}{null} \\[1pt]
\cline{3-3}
&& $e_2^0\cdot f_2^0=\lambda f_1^0$ & \\[1pt]
\cline{3-3}
&& $e_2^0\cdot f_2^0=\lambda f_2^0$ & \\[1pt]
\cline{3-3}
&& $e_2^0\cdot f_2^0=\lambda f_1^0+f_2^0$ & \\[1pt]
\cline{3-3}
&& $e_2^0\cdot f_1^0=\lambda f_1^0+\mu f_2^0,~\mu \neq0$ & \\[1pt]
&& $e_2^0\cdot f_2^0=-\frac{(\lambda-1)\lambda}{\mu}f_1^0+(1-\lambda) f_2^0$ & \\[1pt]
\cline{2-4}
& \multirow{1}{*}{$L_{2|0}^2$} & $e_2^0\cdot f_1^0= f_1^0,~e_2^0\cdot f_2^0= f_2^0$ & null \\[1pt]
\hline
\end{tabular}
}}
\end{description}

\begin{remark}
For the non-graded cases, we recover some of the results obtained in \cite{RE20}.
\end{remark}

\begin{proposition}
Let $(A,L)$ be a Lie-Rinehart superalgebra, with $\dim(A)\leq 2$ and $L$ abelian. Then either the action is trivial, or the anchor is null.
\end{proposition}

\begin{proof}
Since $L$ is abelian, the compatibility condition becomes: $\rho(x)(a)\cdot y=0$ for all $x,y\in L$ and $a\in A$.
\begin{itemize}
\item[$\bullet$] $A=\bf{A_{(1|0)}^1}$. We have already seen that the only compatible pair is the one with the trivial action and the null anchor.

\item[$\bullet$] $A=\bf{A_{(1|1)}^1}$.
f $a=e_1^0$, $\rho(x)(e_1^0)=0$, so the compatibility condition is satisfied. If $a=e_1^1$, we have two cases, depending on the degree of $x$.
\begin{itemize}
\item[$|x|=0$:] it exists $\lambda \in \C$ such that $\rho(x)(e_1^1)=\lambda e_1^1$. For $y\in L$, we obtain $\lambda e_1^1 \cdot y=0$. We then have the following dichotomy: either $\lambda=0$, which means that the anchor is null, or $e_1^1\cdot y=0$, which means that the action is trivial.
\item[$|x|=1$:] it exists $\lambda \in \C$ such that $\rho(x)(e_1^1)=\lambda e_1^0$. For $y\in L$, we obtain $\lambda e_1^0\cdot y=0$, so $\lambda=0$ and the anchor is null.
\end{itemize}

\item[$\bullet$] $A=\bf{A_{(2|0)}^1}$. We have $\rho(x)(e_1^0)=0$ for every $x\in L$. Let $a=e_2^0$. We also have two cases, depending on the degree of $x$:

\begin{itemize}
\item[$|x|=0$:] in this case, there exists $\lambda, \mu\in \C$ such that $\rho(x)(e_2^0)=\lambda e_1^0+\mu e_2^0$. We then have $0=\rho(x)(e_2^0)\cdot y=\lambda y+\mu e_2^0 \cdot y.$ If $\mu=0$, we obtain $\lambda y=0$, so $\lambda=0$ and the anchor is null. If $\mu \neq 0$, we have $ e_2^0\cdot y=-\frac{\lambda}{\mu}y.$

The fourth condition of the Definition \ref{superdef} gives us $\rho(e_2^0\cdot x)(e_2^0)=e_2^0\rho(x)(e_2^0),$ so $-\frac{\lambda^2}{\mu}e_1^0-\lambda e_2^0=\lambda e_2^0.$

We conclude that $\lambda=0$ and the action is trivial.	

\item[$|x|=1$:] $\rho(x)(e_2^0)\in\left(\bf{A_{(2|0)}^1} \right)_1=\left\lbrace 0\right\rbrace$, so the anchor is null. 
\end{itemize}
\item[$\bullet$] $A=\bf{A_{(2|0)}^2}$. The superderivations space is reduced to $\left\lbrace0 \right\rbrace$, so all the anchors are null.
\qedhere
\end{itemize}
\end{proof}

\section{Deformation theory of Lie-Rinehart superalgebras}

In this section, we provide a deformation theory of Lie-Rinehart superalgebras. The following results are strongly inspired by \cite{MM20}, where the authors discussed a deformation theory of Hom-Lie-Rinehart algebras, including Lie-Rinehart algebras. We also use results from \cite{VL15} and \cite{BP89}. One needs a cohomology complex constructed below and that controls deformations.

We could deform four different operations: the multiplication of $A$, the Lie bracket of $L$, the action $A\curvearrowright L$ and the anchor map $\rho:L\longrightarrow \Der(A)$. Here, we restrict ourselves to deforming only the bracket of $L$ and the anchor map $\rho$. It follows that the multiplication of $A$ and the action $A\curvearrowright L$ remain undeformed in the following theory.

All the superalgebras are now $\K$-superalgebras, $\K$ being a characteristic zero field.

\subsection{Super-multiderivations}

Let $A$ be a supercommutative associative $\K$-superalgebra, $M$ an $A$-module and $L$ an $A$-Lie-Rinehart superalgebra with bracket $[\cdot,\cdot]$ and anchor map $\rho$. We recall that both maps are even.

\begin{definition}[Super-multiderivations space]
We define $\Der^n(M,M)$ as the space of multilinear maps $$f:M^{\times n+1}\longrightarrow M$$ such that there exists $\sigma_f:M^{\otimes n}\longrightarrow\Der(A)$ (called symbol map),
\begin{enumerate}
\item For all $i \in \{1, \dotsc, n\}$:
\[
f(x_1,\cdots ,x_i,x_{i+1},\cdots ,x_{n+1})
=-(-1)^{|x_i||x_{i+1}|}f(x_1,\cdots ,x_{i+1},x_i,\cdots ,x_{n+1}),
\]
\item For all $a \in A$:
\begin{align*}
f(x_1,\cdots ,x_n,ax_{n+1})
={ }&{ }(-1)^{|a|(|f|+|x_1|+\cdots +|x_{n}|)}af(x_1,\cdots ,x_{n+1}) \\
&+\sigma_f(x_1,\cdots ,x_{n})(a)(x_{n+1}).
\end{align*}
\end{enumerate}
\end{definition}

\begin{remark}
With this definition, we can check that the bracket $[\cdot,\cdot]$ on $L$ belongs to $\Der^1(L,L)$, with symbol map given by the anchor $\rho$.
\end{remark}
We define
\[\Der^*(M,M)=\bigoplus_{n\geq-1}\Der^n(M,M), \text{ with} \Der^{-1}(M,M)=M.\]
Every space $\Der^n(M,M)$ has a natural $\Z_2$-graduation, given by
\[ |D|=j\in \Z_2\Longleftrightarrow |D(x_1,\cdots ,x_{n+1})|-\sum_i|x_i|=j\text{ mod }2, \text{ for } D\in \Der^n(M,M).  \]
We have
\[ \Der^*(M,M)=\bigoplus_n\left(\Der^n(M,M)\right)= \bigoplus_n\left(\Der^n_0(M,M)\oplus \Der^n_1(M,M)\right).  \]

Next, we provide a bracket on $\Der^*(M,M)$. We adapt the formula of the Nijenhuis-Richardson bracket (\cite{NR67}) to the super case. 
As explained in \cite{VL15}, the space of super-multiderivations $\Der^*(M,M)$ is a $\Z$-graded Lie-algebra, but not a bigraded Lie algebra.
For $f\in\Der^p(M,M)$ and $g\in\Der^q(M,M)$, we define
\begin{align*}
(f\circ g)&(x_1,\cdots ,x_{p+q+1})\\
&=\sum_{\tau\in Sh(q+1,p)}\varepsilon(\tau,x_1,\cdots ,x_{p+q+1})f\left(g(x_{\tau(1)},\cdots ,x_{\tau(q+1)}), x_{\tau(q+2)},\cdots ,x_{\tau(p+q+1)} \right), 
\end{align*}
where $Sh(q+1,p)$ denotes the set of all permutations $\tau$ of $\left\lbrace 0,1,\cdots ,q,\cdots p+q \right\rbrace$ such that \[\tau(0)<\tau(1)<\cdots <\tau(q) \text{ and } \tau (q+1)<\cdots <\tau(p+q).\]	
The sign $\varepsilon(\tau,x_1,\cdots ,x_{p+q+1})$ is implicitly defined by the permutation $\tau$ with respect to the parity of the homogeneous elements $x_1,\cdots ,x_{p+q+1}\in L$. For example, if $\tau=(i,~ i+1)$ is an elementary transposition, then
\[ \varepsilon(\tau,x_1,\cdots ,x_{p+q+1})=sgn(\tau)(-1)^{|x_i||x_{i+1}|}=-(-1)^{|x_i||x_{i+1}|}.  \]
Some computations for $\tau\in S_3$ can be found in \cite{AM10}. An explicit expression of $\varepsilon$ is given in \cite{BP89}.

\begin{proposition}[\text{\cite{VL15}}]	For $f\in\Der^p(M,M)$ and $g\in\Der^q(M,M)$, we define a bracket by
\[
[f,g]=f\circ g-(-1)^{pq} g\circ f,
\]
with symbol map $\sigma_{[f,g]}=\sigma_f\circ g-(-1)^{pq}\sigma_g\circ f +[\sigma_f, \sigma_g]$ and
\[
[\sigma_f, \sigma_g](x_1,\cdots ,x_{p+q})=\sum_{Sh(p,q)}\varepsilon(\tau,x_1,\cdots ,x_{p+q})[\sigma_f(x_{\tau(1)},\cdots ,\tau_{p}),\sigma_g(x_{\tau(p+1)},\cdots ,\tau_{p+q})].  \]
With this bracket, $\Der^*(M,M)$ has a $\Z$-graded Lie algebra structure.
\end{proposition}

\begin{remark}
The reader should be aware that in [VL15], different conventions are adopted, but the result remains unchanged.
\end{remark}

\subsection{Deformation complex}

Let $(A,L)$ be a Lie-Rinehart superalgebra. We construct a deformation complex. 

\begin{proposition}\label{bij}
There is a one-to-one correspondence between Lie-Rinehart superstructures on $(A,L)$ and elements $m\in \Der^1(L,L)$ such that $[m,m]=0$.
\end{proposition}

\begin{proof}
Let $(A,L, [\cdot,\cdot], \rho)$ be a Lie-Rinehart superalgebra. We set $m:=[\cdot,\cdot]$, with symbol map $\sigma_m:=\rho$. We have $m\in Der^1(L,L)$. The super-Jacobi identity for $x_1,x_2,x_3\in L$ is
\[[x_1,[x_2,x_3]]-[[x_1,x_2],x_3]-(-1)^{|x_1||x_2|}[x_2,[x_1,x_3]]=0.  \]

Since $|m|=1$, we have $[m,m](x_1,x_2,x_3)=2 m\circ m(x_1,x_2,x_3).$ Or, by using computations on $S_3$ done in \cite{AM10} for the signs, we get
\begin{align*}
m&\circ m(x_1,x_2,x_3) \\
&=m\left(m(x_1,x_2),x_3 \right) - (-1)^{|x_2||x_3|}m \left(m(x_1,x_3),x_2 \right) + (-1)^{|x_1||x_2|+|x_1||x_3|}m \left(m(x_2,x_3),x_1 \right) \\
&=\left[[x_1,x_2],x_3 \right]+(-1)^{|x_1||x_2|}[x_2,[x_1,x_3]]-[x_1,[x_2,x_3]]\\
&=0.	
\end{align*}
Conversly, let $m\in \Der^1(L,L)$ such that $[m,m]=0$. By setting $[\cdot,\cdot]:=m$ and $\rho:=\sigma_m$, we obtain a Lie-Rinehart superstructure on $(A,L)$. 		
\end{proof}

Hence, we can identify Lie-Rinehart superstructures on $(A,L)$ and the corresponding element $m\in\Der^1(L,L)$. We set
$$C_{def}^n(L,L):=\Der^{n-1}(L,L) \text{ and }
C_{def}^*(L,L):=\bigoplus_{n\geq0}C_{def}^n(L,L),$$
which we endow with an operator 
$$\delta: C_{def}^n(L,L)\longrightarrow C_{def}^{n+1}(L,L),~D\longmapsto [m,D],$$
with  an explicit formula given by 
\begin{align*}
(\delta D)&(x_1,\cdots ,x_{n+1}) \\
={ }&{ }\left( m\circ D-(-1)^{n-1}D\circ m\right)(x_1, \cdots, x_{n+1})\\
={ }&{ }\sum_{\tau\in Sh(n,1)}\varepsilon(\tau,x_1,\cdots, x_{n+1})m\left(D(x_{\tau(1)},\cdots ,x_{\tau(n)}), x_{\tau(n+1)} \right) \\
&-(-1)^{n-1}\sum_{\tau\in Sh(2,n-1)}\varepsilon(\tau,x_1,\cdots, x_{n+1})D\left(m(x_{\tau(1)},x_{\tau(2)}),x_{\tau(3)},\cdots, x_{\tau(n+1)} \right)\\
={ }&{ }\sum_{i=1}^{n+1}\varepsilon_im\left(D(x_1,\cdots ,\hat{x_i},\cdots , x_{n+1}),x_i \right)\\
&-(-1)^{n-1}\sum_{1\leq i<j\leq n} \varepsilon_i^jD\left(m(x_i,x_j), x_1,\cdots ,\hat{x_i},\cdots ,\hat{x_j},\cdots ,x_{n+1} \right), 
\end{align*}
where $\varepsilon_i$ and $\varepsilon_i^j$ denote the signs associated to the permutations with respect to the parity of the homogeneous elements $x_1,\cdots ,x_{n+1}\in L$ and $D\in\Der^*(L,L)$.

\begin{proposition}
The operator $\delta$ is a differential, $i.e.$ $\delta^2=0$.
\end{proposition}

\begin{proof} Let $D\in C_{def}^n(L,L)$. We have
\[
\delta^2(D)
=[m,[m,D]]
=[[m,m],D]+(-1)^{|m||m|}[m,[m,D]]
=0-[m,[m,D]].
\]
Then $[m,[m,D]]=-[m,[m,D]]$, so $[m,[m,D]]=0$.	
\end{proof}

This enables us to define a cohomology complex, which will be used next in the deformation theory of Lie-Rinehart superalgebras. Therefore we call it \textbf{deformation cohomology}. For $p,q \in \Z$, we have the usual definitions of $p$-cocycles and $q$-coboundaries, denoted respectively by $Z_{def}^p(L)$ and $B_{def}^q(L)$. Finally we set $H^p_{def}(L):=Z_{def}^p(L)/B_{def}^p(L)$ to be the $p$-th cohomology group. We have
\begin{align*}
Z^1_{def}(L)
&=\ker(\delta^1) \\
&=\left\{ D\in \Der^0(L),~D\left( \left[ x,y\right] \right) =\left[ D(x),y \right] -(-1)^{|x||y|}\left[ D(y),x \right]~\forall x,y\in L \right\},
\end{align*}
and $Z^2_{def}(L) = \ker(\delta^2)$ is the set consisting of $D \in \Der^1(L)$ which satisfies the condition:
\begin{align*}
\left[ D(x,y),z\right] { }&{ } -(-1)^{|y||z|}\left[ D(x,z),y \right] +(-1)^{|x||y|+|x||z|}\left[D(y,z),x \right] \\
&{ }= -D\left([x,y],z \right) + (-1)^{|y||z|}D\left([x,z],y \right) + (-1)^{|x||y|+|x||z|}D\left( [y,z],x \right),
\end{align*}
for all $x,y,z\in L$.

\subsection{Formal deformations}

In this section, we discuss deformation theory of Lie-Rinehart superalgebras and show that the deformations are controlled by the cohomology defined above. Notice that in the sequel, we aim to deform the Lie bracket and the anchor while we keep fixed the multiplication of the associative superalgebra and its action. We denote by $\K[[t]]$ (resp. L[[t]]) the formal power series ring in $t$ with coefficients in $\K$ (resp. the formal space in $t$ with coefficients in the vector superspace $L$).

\begin{definition}
Let $(A,L,[\cdot,\cdot],\rho)$ be a Lie-Rinehart superalgebra over a field $\K$ of characteristic zero, and let $m\in \Der^1(L,L)$ be the corresponding element obtained by Proposition~\ref{bij}. A deformation of the Lie-Rinehart superalgebra is given by a $\K[[t]]$-bilinear map
$$m_t:L\times L\longrightarrow L[[t]],~m_t(x,y)=\sum_{i\geq 0}t^im_i(x,y),$$
such that $m_0=m \text{ and } m_i\in\Der^1(L,L)$ with symbol map denoted by $\sigma_{m_i}$ for $i \ge 1$, satisfying $[m_t,m_t]=0$, the bracket being the $\Z$-graded bracket on $\Der^*(L[[t]],L[[t]])$.
\end{definition}

\begin{remark}
The map $m_t$ defined on $L\times L$ can be extended to a map on $L[[t]]\times L[[t]]$ using the $\K[[t]]$-bilinearity.
\end{remark}

We check that $m_t$ is a 1-degree super-multiderivation of $L[[t]]$, with symbol map given by $\sigma_{m_t}=\sum_it^i\sigma_{m_i}$. As a consequence, $m_t$ gives rise to a Lie-Rinehart superstructure on $(A[[t]], L[[t]])$, with bracket $[\cdot,\cdot]_t:=m_t$ and anchor $\rho_t:=\sigma_{m_t}$.

\begin{remark}
The first non-zero element $m_i$ of the deformation is called the \textbf{infinitesimal} of the deformation.
\end{remark}

Since $m_t$ satisfies $[m_t,m_t]=0$, we have
\begin{equation}\label{eqdef}
m_t(x_1,m_t(x_2,x_3))=m_t(m_t(x_1,x_2),x_3)+(-1)^{|x_1||x_2|}m_t(x_2,m_t(x_1,x_3)).  
\end{equation}
This equation is called deformation equation and is equivalent to an infinite system by identifying the coefficients of $t$.

\begin{theorem}
Let $m_t$ be a deformation of a Lie-Rinehart superalgebra $(A,L)$. Then the infinitesimal of the deformation $m_1$ is a 2-cocycle with respect to the deformation cohomology.
\end{theorem}

\begin{proof}
By taking the coefficients of $t$ in Equation (\ref{eqdef}), we obtain
\begin{align*}
&~~~~ m_1(x_1,m(x_2,x_3))-m_1(m(x_1,x_2),x_3)-(-1)^{|x_1||x_2|}m_1(x_2,m(x_1,x_3))\\
&+m(x_1,m_1(x_2,x_3))-m(m_1(x_1,x_2),x_3)-(-1)^{|x_1||x_2|}m(x_2,m_1(x_1,x_3))=0.
\end{align*}
Using \cite{AM10} for the signs, we get 
\begin{align*}
[m,m_1](x_1,x_2,x_3)
={ }&{ }m(m_1(x_1,x_2),x_3) - (-1)^{|x_2||x_3|}m(m_1(x_1,x_3),x_2) \\
&+(-1)^{|x_1||x_2|+|x_1||x_3|}m(m_1(x_2,x_3),x_1) +m_1(m(x_1,x_2),x_3) \\
&-(-1)^{|x_2||x_3|}m_1(m(x_2,x_3),x_1)+(-1)^{|x_1||x_2|+|x_1||x_3|}m_1(m(x_2,x_3),x_1)\\
={ }&{ }m(m_1(x_1,x_2),x_3) + (-1)^{|x_1||x_2|}m(x_2,m_1(x_1,x_3)) \\
&-m(,x_1,m_1(x_2,x_3))+m_1(m(x_1,x_2),x_3) \\
&+(-1)^{|x_1||x_2|}m_1(x_2,m(x_1,x_3))-m_1(x_1,m(x_2,x_3))\\
={ }&{ }0.
\qedhere
\end{align*}
\end{proof}

\subsection{Equivalence of deformations}
Let $(A,L,[\cdot,\cdot],\rho)$ be a Lie-Rinehart superalgebra and $m$ the associated element of $\Der(L,L)$. Let $m_t$ and $m_t'$ be two deformations of $m$.

\begin{definition}
We say that $m_t$ and $m_t'$ are \textbf{equivalent} if there exists an even formal automorphism $\Phi_t$ of $L[[t]]$, that can be written $ \Phi_t=\text{id}+\sum_{i\geq 1}t^i\phi_i, \text{ with }\phi_i:L\longrightarrow L$ even $\K$-linear maps, such that $\Phi_t\circ m_t'(x,y)=m_t(\Phi_t(x),\Phi_t(y)). $ We write $m_t\sim m_t'$.
\end{definition}

\begin{definition}
A deformation is said to be \textbf{trivial} if it is equivalent to the deformation given by $m_t^0=\sum t^im_i^0$, with $m_0^0=m$ and $m_i^0=0$ for $i\geq 1$.
\end{definition}

We recall that we have a short exact sequence $$ 0\longrightarrow B^2_{def}(L)\longrightarrow Z^2_{def}(L)\overset{\pi}{\longrightarrow} H^2_{def}(L)\longrightarrow0.    $$
We denote $\overline{\mu}=\pi(\mu)$ for $\mu\in Z^2_{def}(L)$.

\begin{theorem}
For a deformation $m_t$ of $m$, the cohomology class of the infinitesimal element $m_1$ is determined by the equivalence class of $m_t$.
\end{theorem}

\begin{remark}
In other words, we have $m_t\sim m_t'\Longrightarrow \overline{m_1}=\overline{m_1'}$.
\end{remark}

\begin{proof}
Let $m_t$ and $m_t'$ be two equivalent deformations of $m$ and $\Phi_t$ the associated formal automorphism. By definition, we have $\Phi_t\circ m_t'(x,y)=m_t(\Phi_t(x),\Phi_t(y)), $ that can be rewritten $ \sum_{k,i}t^{k+1}\phi_k(m_i'(x,y))=\sum_{j,p,q}t^{j+p+q}m_j\left( \phi_p(x),\phi_q(y)\right).$ By identifying coefficients of $t$, we obtain	
\[ m_1(x,y)-m_1'(x,y)=\phi_1(m(x,y))-m(\phi_1(x),y)-m(x,\phi_1(y)).  \]

Since $\delta(\phi_1)=m(\phi_1(x),y)+m(x,\phi_1(y))-\phi_1(m(x,y)),  $ we have $m_1'-m_1=\delta(\phi_1)$. It follows that $m_1'=m_1+\delta(\phi_1)$, so $\overline{m_1}=\overline{m_1'}\in H^2_{def}(L)$.
\end{proof}

\begin{definition}
A Lie-Rinehart superalgebra is said to be \textbf{rigid} if every deformation is trivial.
\end{definition}

\begin{theorem}
Any non-trivial deformation of $m\in Der^1(L,L)$ is equivalent to a deformation whose infinitesimal is not a coboundary.
\end{theorem}

\begin{remark}
That can be reformulated as if all the elements $m_i$ are coboundaries, then $m_t\sim m_t^0$.
\end{remark}

\begin{proof}
Suppose $m_1$ is a coboundary: $~\exists \phi\in C^1_{def}=\Der^0(L,L)$ such that $m_1=\delta(\phi)$. We show that $\overline{m_1}=0$. We set $\Phi_t=\text{id}+t\phi$ and define $m_t':=\Phi_t\circ m_t \circ \Phi_t^{-1}$. Then we have $m_t\sim m_t'$, that is 
\[ \sum_j t^jm_j(\left( \Phi_t(x),\Phi_t(y) \right) = \Phi_t\left( \sum_i t^im_i'(x,y) \right), \] which is equivalent to
\[\sum_{j,k,l}t^{j+k+l}m_j\left( \phi_k(x), \phi_l(y) \right)=\sum_{i,p}t^{i+p}\phi_p(m_i'(x,y)).   \]
By identifying the coefficients of $t$, we get
\[
m_1'(x,y)-m_1(x,y)=\phi(m(x,y))-m'(\phi(x),y)-m'(x,\phi(y))
=-\delta(\phi).
\]
Then, we have $m_1'-m_1=-\delta(\phi)=-m_1$, so $m_1'=0$. By repeating the argument, we show that if $m_i\in B^2$, then $\overline{m_i}=0$.	
\end{proof}

\begin{corollary}
If $H^2_{def}(L)=0$, any deformation is equivalent to a trivial deformation.
\end{corollary}

\begin{proof}
If $H^2_{def}(L)=0$, the infinitesimal is a coboundary. According to the theorem, the deformation is equivalent to a trivial deformation.
\end{proof}

\subsection{Obstructions}

Let $(A,L,[\cdot,\cdot],\rho)$ be a Lie-Rinehart superalgebra, $m$ the associated element of $\Der(L,L)$ and $N\in \N,$ $N\geq1$. We say that $m_t$ is a deformation of $m$ of order $N$ or $N$-order deformation if
\[m_t=\sum_{k=0}^{N} t^km_k,~~m_k\in\Der^1(L,L)~~\text{and }[m_t,m_t]=0. \]
Here we aim to extend a $N$-order deformation $m_t$ to a $(N+1)$-order deformation, $i.e.$ find $m_{N+1}\in\Der^1(L,L)$ such that $m_t'=m_t+t^{N+1}m_{N+1}$ is a deformation of $m$.

The condition on $m_{N+1}$ is expressed by the following deformation equation 
\[\delta m_{N+1}(a,b,c)=\sum_{\underset{i,j>0}{i+j=N}}m_i(a,m_j(b,c))-m_i(m_j(a,b),c)-(-1)^{|a||b|}m_i(b,m_j(a,c)).   \]
\begin{definition}
We set, for $a,b,c\in L$,
\[\theta_{N}(a,b,c)=\sum_{\underset{i,j>0}{i+j=N}}m_i(a,m_j(b,c))-m_i(m_j(a,b),c)-(-1)^{|a||b|}m_i(b,m_j(a,c)).   \]
It's immediate that $\theta_{N}\in C_{def}^3(L,L)=\Der^2(L,L)$. The map $\theta_{N}$ is called the obstruction cochain of the $N$-order deformation $m_t$.
\end{definition}

\begin{lemma}
\[\theta_{N}=-\frac{1}{2}\sum_{\underset{i,j>0}{i+j=N}}[m_i,m_j].\]
\end{lemma}

\begin{corollary}
The map $\theta_{N}$ is a 3-cocycle.
\end{corollary}

\begin{proof}
Using graded Jacobi, we have
\[
\delta(\theta_{N})=[m,\theta_{N}]=-\frac{1}{2}\sum_{\underset{i,j>0}{i+j=N}}\left[ m,[m_i,m_j] \right]=-\frac{1}{2}\sum_{\underset{i,j>0}{i+j=N}}\left[ [m,m_i],m_j \right]+\frac{1}{2}\sum_{\underset{i,j>0}{i+j=N}}\left[ m_i,[m,m_j] \right].
\]
Since $[m,\theta_{N}]=0$ if and only if $\sum_{\underset{i,j>0}{i+j=N}}\left[ [m,m_i],m_j \right]=\sum_{\underset{i,j>0}{i+j=N}}\left[ m_i,[m,m_j] \right]$, it follows that $\sum_{\underset{i,j>0}{i+j=N}}\left[ m_j,[m,m_i] \right]=\sum_{\underset{i,j>0}{i+j=N}}\left[ m_i,[m,m_j] \right]$.
\end{proof}\begin{theorem}
Let $m_t$ be a $N$-order deformation of $m$. Then $m_t$ extends to a $(N+1)$-order deformation if and only if $\theta_{N}$ is a 3-coboundary.
\end{theorem}

\begin{proof}
\begin{enumerate}
\item[$(\Rightarrow)$] Suppose that $m_t'$ is a $(N+1)$-order deformation of $m$. Then $m_t'$ satisfies the graded Jacobi identity, that is for $a,b,c\in L$:

\[ m_t'(a,m_t'(b,c))-m_t'(m_t'(a,b),c)-(-1)^{|a||b|}m_t'(b,m_t'(a,c))=0. \]
By expanding and collecting the coefficients of $t^{N+1}$, we have		

\[\sum_{\underset{i,j\geq0}{i+j=N+1}}m_i(a,m_j(b,c))-m_i(m_j(a,b),c)-(-1)^{|a||b|}m_i(b,m_j(a,c))=0,   \] which is equivalent to	$-[m,m_{N+1}]+\theta_{N}(a,b,c)=0.$	As a consequence, 
\[\theta_{N}(a,b,c)=\delta(m_{N+1}).\]

\item[$(\Leftarrow)$] 	Suppose $\theta_{N}$ is a coboundary: it exists $\varphi\in C_{def}^2(L)$ such that $\theta_{N}=\delta \varphi=[m,\varphi]$. We need to show that $m_t'=m_t+t^{N+1}\varphi$ is a $(N+1)$-order deformation of $m$.
We write
\[ \sum_{\underset{i,j>0}{i+j=N+1}}m_i(a,m_j(b,c))-m_i(m_j(a,b),c)-(-1)^{|a||b|}m_i(b,m_j(a,c))=[m,\varphi],\] 
and it follows that
\[ \sum_{\underset{i,j\geq0}{i+j=N+1}}m_i(a,m_j(b,c))-m_i(m_j(a,b),c)-(-1)^{|a||b|}m_i(b,m_j(a,c))=[m,\varphi]=0. \]	
With this equality and with the fact that $m_t$ already satisfies the super-Jacobi identity, we deduce that $m_t'$ is a $(N+1)$-order deformation of $m$.
\qedhere
\end{enumerate}
\end{proof}

\begin{corollary}
If $H_{def}^3(L)=0$, any $N$-order deformation extends to a $(N+1)$-order deformation.
\end{corollary}

\subsection{Rigid $(1|1,1|1)$-type Lie-Rinehart superalgebra}

In order to provide an example of a rigid Lie-Rinehart superalgebra, consider the pair $\left( \mathbf{A_{1|1}^1}, \mathbf{L_{1|1}^1} \right)$ (product on $\bf{A_{1|1}^1}:$ $ e_1^1e_1^1=0$; bracket on $\bf{L_{1|1}^1}:$ $[f_1^1,f_1^1]=f_1^0$), endowed with the null anchor and the action given by $e_1^1\cdot f_1^1=\lambda f_1^0,~\lambda\in \C$.
Using Proposition \ref{bij}, the bracket $[\cdot,\cdot]$ corresponds to $m\in\Der^1(L)$, with symbol map $\sigma_m=\rho$. We will show that $(A,L)$ endowed with this Lie-Rinehart superstructure is rigid. 

\begin{lemma}\label{fidele}
Let $(A,L)$ being endowed with the structure above and let $D\in\Der^n(L)$. If $D=0$, then $\sigma_D(X)(a)=0$ for all $X\in L^{\times n}$ and all $a\in A$. 
\end{lemma}

Now we aim to compute explicitly elements of $Z^2_{def}(L)$. We set the general form of a 2-cochain $(D,\sigma_D)$:
\[
\begin{cases}
D(f_1^0,f_1^0)=0,\\
D(f_1^0,f_1^1)=\gamma_0f_1^0+\gamma_1f_1^1,\\
D(f_1^1,f_1^1)=\theta_0f_1^0+\theta_1f_1^1,
\end{cases}
\qquad
\begin{cases}
\sigma_D(f_1^0)(e_1^0)=\sigma_D(f_1^1)(e_1^0)=0,\\
\sigma_D(f_1^0)(e_1^1)=p_0e_1^0+p_1e_1^1,\\
\sigma_D(f_1^1)(e_1^1)=q_0e_1^0+q_1e_1^1.
\end{cases}
\]	
All the parameters belong to $\C$. The following result gives us conditions on these parameters for $(D,\sigma_D)$ to belong to $Z^2_{def}(L)$.

\begin{lemma}\label{Z2}
Let $(A,L)$ be a Lie-Rinehart superalgebra endowed with the above structure and let $(D,\sigma_D)\in Z^2_{def}(L)$. Then
\[
\begin{cases}
D(f_1^0,f_1^0)=0,\\
D(f_1^0,f_1^1)=\gamma f_1^0,\\
D(f_1^1,f_1^1)=\theta f_1^0-\gamma f_1^1,
\end{cases}
\qquad
\begin{cases}
\sigma_D(f_1^0)(e_1^0)=\sigma_D(f_1^1)(e_1^0)=0,\\
\sigma_D(f_1^0)(e_1^1)=0,\\
\sigma_D(f_1^1)(e_1^1)=q_0e_1^0+q_1e_1^1.
\end{cases}
\]
\end{lemma}

\begin{proof}
Suppose $(D,\sigma_D)\in Z^2_{def}(L)$. Then $\delta(D)=[m,D]=0.$ Evaluating this equation on the basis elements $f_1^0,f_1^1$ of $L$, we find 
\[
[m,D](f_1^1,f_1^0,f_1^1)=-2\gamma_1 f_1^0
\quad \textup{and} \quad
[m,D](f_1^1,f_1^1,f_1^1)=(\gamma_0+\theta_1) f_1^0+\gamma_1 f_1^1,
\]
all other possible combinations being zero. Now if $u=u_0f_1^0+u_1f_1^1$, $v=v_0f_1^0+v_1f_1^1$ and $w=w_0f_1^0+w_1f_1^1$, we have
\[ [m,D](u,v,w)=-2u_1v_0w_1\gamma_1f_1^0+u_1v_1w_1\left((\gamma_0+\theta_1)f_1^0+\gamma_1f_1^1 \right)=0.  \]
Setting $t_0:=-2u_1v_0w_1$ and $t_1:=u_1v_1w_1$, we get
\begin{center}
$\begin{cases}
t_0\gamma_1+t_1(\gamma_0+\theta_1)&=0\\
t_1\gamma_1&=0.
\end{cases}$
\end{center}
As a consequence, $\gamma_1=0$ and $\gamma_0=-\theta_1$. One gets the expression in the Lemma by setting  $\gamma_0=:\gamma$ and $\theta_0=:\theta$.

Then, we know by Lemma \ref{fidele} that $\sigma_{[m,D]}(z)(a)=0$ for any $z\in L $ and $a\in A$. In this case $\sigma_{[m,D]}=\rho\circ D+\sigma_D\circ m+[\rho,\sigma_D]=\sigma_D\circ m,$ because $\rho=0$. Therefore, we have
\begin{center}
$\begin{cases}
\sigma_{[m,D]}(f_1^0,f_1^0)=0\\
\sigma_{[m,D]}(f_1^0,f_1^1)=0\\
\sigma_{[m,D]}(f_1^1,f_1^1)(e_1^1)=\sigma_D(f_1^0)(e_1^1)=p_0e_1^0+p_1e_1^1.\\
\end{cases}$
\end{center}
Since $\sigma_{[m,D]}=0,$ we obtain $p_0=p_1=0$.
\end{proof}

For $\Delta\in C^1_{def}(L)$, we have $\sigma_{\Delta}=0$. We aim to describe $B^2_{def}(L)$. We write for $\lambda_0, \lambda_1,\mu_0,\mu_1\in\C$:

\begin{center}
$\begin{cases}
\Delta(f_1^0)=\lambda_0f_1^0+\lambda_1f_1^1,\\
\Delta(f_1^1)=\mu_0f_1^0+\mu_1f_1^1.\\
\end{cases}$
\end{center} 

Then we have:
\begin{align*}
[m,\Delta](f_1^0,f_1^0)&=0,\\
[m,\Delta](f_1^0,f_1^1)&=\lambda_2x,\\
[m,\Delta](f_1^1,f_1^1)&=(2\mu_1-\lambda_0)f_1^0-\lambda_1f_1^1.\\
\end{align*}

\begin{proposition}
If $(A,L)$ is endowed with the above structure, then $B^2_{def}(L)=Z^2_{def}(L)$, that is, $H^2_{def}(L)=\{0\}$. Thus, $(A,L)$ is rigid.
\end{proposition}

\begin{proof}
Let $D\in Z^2_{def}(L)$, given by the Lemma \ref{Z2}. We need to find $\Delta\in C^1_{def}(L)$ such that 
\[[m,\Delta]=D~\text{ and }~ \sigma_{[m,\Delta]}=\sigma_D.\]
If we set
\begin{center}
$\begin{cases}
\Delta(f_1^0)=\lambda_0f_1^0+\gamma f_1^1\\
\Delta(f_1^1)=\mu_0f_1^0+\frac{\theta+\lambda_0}{2}f_1^1,\\
\end{cases}$
\end{center} 
for arbitrary $\lambda_0$ and $\mu_0$, we have $[m,\Delta]=D$. Because $\sigma_{\Delta}=0$, we have $\sigma_{[m,\Delta]}=\rho\circ\Delta$. But both those terms are zero, so this last equation is verified.
\end{proof}


\section{Appendix}
\subsection{Supercommutative associative superalgebras}\label{appendixasso}
We list supercommutative associative superalgebras with unit. We denote basis elements of $A_0$ by $e_i^0$ and those of $A_1$ by $e_j^1$. The unit is $e_1^0$. 

\begin{itemize}
\item[$\bullet$] The purely odd superalgebras $\bf{A_{0|p}}$ always have a zero product.

\item[$\bullet$] $\dim A=(1|0)$: there is only one unital supercommutative associative superalgebra $\bf{A_{1|0}^1}$, with product $e_1^0e_1^0=e_1^0.$

\item[$\bullet$] $\dim A=(1|1)$: there is only one unital supercommutative associative superalgebra $\bf{A_{1|1}^1}$ with product $e_1^1e_1^1=0$.

\item[$\bullet$] $\dim A=(1|p),~p\geq 2$: we have $(A_1)^2=\left\lbrace 0\right\rbrace $ (\cite{ACZ09}).	

\item[$\bullet$] $\dim A=(2|0)$: there are two pairwise non-isomorphic unital supercommutative associative superalgebras:
\begin{itemize}
\item [$\bf{A_{2|0}^1}:$] every non-trivial product is zero;
\item [$\bf{A_{2|0}^2}:$] $e_2^0e_2^0=e_2^0$.
\end{itemize}

\item[$\bullet$] $\dim A=(2|1)$: there are two pairwise non-isomorphic unital supercommutative associative superalgebras:
\begin{itemize}
\item [$\bf{A_{2|1}^1}:$] every non-trivial product is zero;
\item [$\bf{A_{2|1}^2}:$] $e_2^0e_2^0=e_2^0,~e_2^0e_1^1=e_1^1$.
\end{itemize}

\item[$\bullet$] $\dim A=(2|2)$: there are five pairwise non-isomorphic unital supercommutative associative superalgebras:
\begin{itemize}
\item [$\bf{A_{2|2}^1}:$] $e_2^0e_2^0=e_2^0,~e_2^0e_1^1=e_1^1$;
\item [$\bf{A_{2|2}^2}:$] $e_2^0e_2^0=e_2^0$;
\item [$\bf{A_{2|2}^3}:$] $e_2^0e_1^1=e_2^1$;
\item [$\bf{A_{2|2}^4}:$] every non-trivial product is zero;
\item [$\bf{A_{2|2}^5}:$] $e_1^1e_2^1=e_2^0$.
\end{itemize}

\item[$\bullet$] $\dim A=(3|0)$: there are four pairwise non-isomorphic unital supercommutative associative superalgebras:
\begin{itemize}
\item [$\bf{A_{3|0}^1}:$] $e_2^0e_2^0=e_2^0,~e_2^0e_3^0=e_3^0,~e_3^0e_3^0=e_3^0$;
\item [$\bf{A_{3|0}^2}:$] $e_2^0e_2^0=e_2^0,~e_2^0e_3^0=e_3^0$;
\item [$\bf{A_{3|0}^3}:$] $e_2^0e_2^0=e_2^0$;
\item [$\bf{A_{3|0}^4}:$] every non-trivial product is zero;
\end{itemize}

\item[$\bullet$] $\dim A=(3|1)$: there are five pairwise non-isomorphic unital supercommutative associative superalgebras:
\begin{itemize}
\item [$\bf{A_{3|1}^1}:$] $e_2^0e_2^0=e_2^0,~e_3^0e_3^0=e_3^0$;
\item [$\bf{A_{3|1}^2}:$] $e_2^0e_2^0=e_2^0,~e_2^0e_3^0=e_3^0$;
\item [$\bf{A_{3|1}^3}:$] $e_2^0e_2^0=e_2^0$;
\item [$\bf{A_{3|1}^4}:$] $e_2^0e_2^0=e_3^0$;
\item [$\bf{A_{3|1}^5}:$] every non-trivial product is zero;
\end{itemize}

\item[$\bullet$] $\dim A=(4|0)$: there are nine pairwise non-isomorphic unital supercommutative associative superalgebras:
\begin{itemize}
\item [$\bf{A_{4|0}^1}:$] $e_2^0e_2^0=e_2^0,~e_3^0e_3^0=e_3^0,~e_4^0e_4^0=e_4^0$;
\item [$\bf{A_{4|0}^2}:$] $e_2^0e_2^0=e_2^0,~e_3^0e_3^0=e_3^0$;
\item [$\bf{A_{4|0}^3}:$] $e_2^0e_2^0=e_2^0~e_2^0e_3^0=e_3^0$;
\item [$\bf{A_{4|0}^4}:$] $e_2^0e_2^0=e_2^0,~e_2^0e_3^0=e_3^0,~e_2^0e_4^0=e_4^0,~e_3^0e_3^0=e_4^0$;
\item [$\bf{A_{4|0}^5}:$] $e_2^0e_2^0=e_3^0,~e_2^0e_3^0=e_4^0$;
\item [$\bf{A_{4|0}^6}:$] $e_2^0e_2^0=e_2^0,~e_2^0e_3^0=e_3^0,~e_2^0e_4^0=e_4^0$;
\item [$\bf{A_{4|0}^7}:$] $e_2^0e_3^0=e_4^0$;
\item [$\bf{A_{4|0}^8}:$] $e_3^0e_3^0=e_3^0$;
\item [$\bf{A_{4|0}^9}:$] every non-trivial product is zero.
\end{itemize}
\end{itemize}

\subsection{Lie-Rinehart superalgebras with $\dim(L)>2$}\label{appendix2}

\begin{table}
\centering
\begin{tabular}{|c|c|l|l|}
\hline
$\textbf{A}$ & $\textbf{L}$ & $\textbf{Action}$ & $\textbf{Anchor}$ \\
\hline
\multirow{8}{*}{$A_{1|1}^1$} & \multirow{3}{*}{$L_{1|2}^1$} & trivial & $\rho(f_1^0)(e_1^1)=\lambda e_1^1$\\
\cline{3-4}
& & $e_1^1\cdot f_1^0=\mu f_1^1$ & $\rho(f_1^0)(e_1^1)=e_1^1$ \\
\cline{3-4}
&& $e_1^1\cdot f_1^0=\mu f_2^1$ & $\rho(f_1^0)(e_1^1)=pe_1^1$ \\
\cline{2-4}
& \multirow{2}{*}{$L_{1|2}^2$}& trivial &$\rho(f_1^0)(e_1^1)=\lambda e_1^1$\\
\cline{3-4}
& & $e_1^1\cdot f_1^0=\mu f_1^1$ & null \\
\cline{2-4}
& \multirow{2}{*}{$L_{1|2}^3$}& trivial & $\rho(f_1^0)(e_1^1)=\lambda e_1^1$ \\
\cline{3-4}
\multirow{5}{*}{$A_{1|1}^1$}& & $e_1^1\cdot f_1^0=\mu f_1^1$ & $\rho(f_1^0)(e_1^1)= e_1^1$ \\
\cline{2-4}
& \multirow{3}{*}{$L_{1|2}^4$}& triviale & $\rho(f_1^0)(e_1^1)= \lambda e_1^1$ \\
\cline{3-4}
& & $e_1^1\cdot f_1^0=\mu (f_1^1-if_2^1)$ & $\rho(f_1^0)(e_1^1)= (p-i) e_1^1$ \\
\cline{3-4}
& & $e_1^1\cdot f_1^0=\mu (f_1^1+if_2^1)$ & $\rho(f_1^0)(e_1^1)= (p+i) e_1^1$ \\
\cline{2-4}
& \multirow{1}{*}{\smash{$L_{1|2}^5$}}& \smash{$e_1^1\cdot f_1^1=\mu f_1^0,~e_1^1f_2^1=\gamma f_1^0$} & null\\
\cline{2-4}
& \multirow{3}{*}{$L_{1|2}^6$}& trivial & $\rho(f_1^1)(e_1^1)= \lambda e_1^0$\\
\cline{3-4}
&&$e_1^1\cdot f_1^1=\lambda f_1^0,~e_1^1\cdot f_2^1=\mu f_1^0$&null\\
\cline{3-4}
&&$e_1^1\cdot f_1^0=\lambda f_1^1+\mu f_2^1$&null\\
\hline
\end{tabular}
\caption{$(1|1, 1|2)$-type, $(\lambda, \mu,\gamma \in\C)$}
\end{table}

\begin{table}
\centering
\begin{tabular}{|c|c|l|l|}
\hline
$\textbf{A}$ & $\textbf{L}$ & $\textbf{Action}$ & $\textbf{Anchor}$ \\
\hline
\multirow{15}{*}{$A_{1|1}^1$} & \multirow{4}{*}{$L_{1|3}^1$} & trivial & $\rho(f_1^0)(e_1^1)=\lambda e_1^1$\\
\cline{3-4}
& & $e_1^1\cdot f_1^0=\mu f_1^1$ & $\rho(f_1^0)(e_1^1)=e_1^1$ \\
\cline{3-4}
& & $e_1^1\cdot f_1^0=\mu f_2^1$ & $\rho(f_1^0)(e_1^1)=pe_1^1$ \\
\cline{3-4}
& & $e_1^1\cdot f_1^0=\mu f_3^1$ & $\rho(f_1^0)(e_1^1)=qe_1^1$ \\
\cline{2-4}
& & trivial & $\rho(f_1^0)(e_1^1)=\lambda e_1^1$\\
\cline{3-4}
& $L_{1|3}^2$& $e_1^1\cdot f_1^0=\mu f_2^1$ & null\\
\cline{3-4}
& & $e_1^1\cdot f_1^0=\mu f_1^1$ & $\rho(f_1^0)(e_1^1)=e_1^1$\\
\cline{2-4}
& & trivial & $\rho(f_1^0)(e_1^1)=\lambda e_1^1$ \\
\cline{3-4}
& $L_{1|3}^3$& $e_1^1\cdot f_1^0=\mu f_1^1$ & $\rho(f_1^0)(e_1^1)= pe_1^1$ \\
\cline{3-4}
& & $e_1^1\cdot f_1^0=\mu f_2^1$ & $\rho(f_1^0)(e_1^1)=e_1^1$ \\
\cline{2-4}
& \multirow{1}{*}{\smash{$L_{1|3}^4$}}& trivial & \smash{$\rho(f_1^0)(e_1^1)= \lambda e_1^1$} \\
\cline{3-4}
&\multirow{3}{*}{$L_{1|3}^4$} & $e_1^1\cdot f_1^0=\mu f_1^1$ & $\rho(f_1^0)(e_1^1)= pe_1^1$ \\
\cline{3-4}
& & $e_1^1\cdot f_1^0=\gamma (f_2^1-if_3^1)$ & $\rho(f_1^0)(e_1^1)= (q-i)e_1^1$\\
\cline{3-4}
&& $e_1^1\cdot f_1^0=\gamma (f_2^1+if_3^1)$ & $\rho(f_1^0)(e_1^1)= (q+i)e_1^1$\\
\cline{2-4}
& \multirow{2}{*}{$L_{1|3}^5$}& trivial & $\rho(f_1^0)(e_1^1)= \lambda e_1^1$\\
\cline{3-4}
& & $e_1^1\cdot f_1^0=\mu f_1^1$ & null\\
\cline{2-4}
& \multirow{2}{*}{$L_{1|3}^6$}& trivial & $\rho(f_1^0)(e_1^1)= \lambda e_1^1$\\
\cline{3-4}
& & $e_1^1\cdot f_1^0=\mu f_1^1$ & $\rho(f_1^0)(e_1^1)= e_1^1$\\
\cline{2-4}
& $L_{1|3}^7$&$e_1^1\cdot f_1^1=\mu f_1^0,~e_1^1\cdot f_2^1=\gamma f_1^0,~e_1^1\cdot f_3^1=\lambda f_1^0$&null\\
\cline{2-4}
& \multirow{4}{*}{$L_{1|3}^8$}& trivial & $\rho(f_1^0)(e_1^1)= \lambda e_1^1$\\
\cline{3-4}
& & $e_1^1\cdot f_1^1=\lambda f_1^0,~e_1^1\cdot f_2^1=\mu f_1^0$ & \multirow{2}{*}{null}\\&&$e_1^1\cdot f_3^1=\gamma f_1^0$&\\
\cline{3-4}
&&$e_1^1\cdot f_1^0=\lambda f_1^1+\mu f_2^1+\gamma f_3^1$&null\\
\hline
\end{tabular}
\caption{$(1|1, 1|3)$-type, $(\lambda, \mu,\gamma \in\C)$}
\end{table}

\begin{table}
\centering
\begin{tabular}{|c|c|l|l|}
\hline
$\textbf{A}$ & $\textbf{L}$ & $\textbf{Action}$ & $\textbf{Anchor}$ \\
\hline
\multirow{13}{*}{$A_{1|1}^1$} & \multirow{2}{*}{$L_{2|1}^1$} & trivial & $\rho(f_1^0)(e_1^1)=\lambda e_1^1$ \\
\cline{3-4}
& & $e_1^1\cdot f_1^1=\lambda f_1^0+\mu f_2^0$ & null \\
\cline{2-4}
& \multirow{2}{*}{$L_{2|1}^2$} & trivial & $\rho(f_1^0)(e_1^1)=\lambda e_1^1,~\rho(f_2^0)(e_1^1)=\mu e_1^1$ \\
\cline{3-4}
& & $e_1^1\cdot f_1^0=\lambda f_1^1,~e_1^1\cdot f_2^0=\mu f_1^1$ & $\rho(f_1^0)(e_1^1)= e_1^1,~\rho(f_2^0)(e_1^1)= -e_1^1$ \\
\cline{2-4}
& $L_{2|1}^3$ & trivial & $\rho(f_1^0)(e_1^1)=\lambda e_1^1$ \\
\cline{2-4}
& \multirow{2}{*}{$L_{2|1}^4$} & trivial & $\rho(f_1^0)(e_1^1)=\lambda e_1^1$ \\
\cline{3-4}
& & $e_1^1\cdot f_1^0=\mu f_1^1$ & $\rho(f_1^0)(e_1^1)=pe_1^1$ \\
\cline{2-4}
& & trivial & $\rho(f_1^0)(e_1^1)=\lambda e_1^1$ \\
\cline{3-4}
& $L_{2|1}^5$ & $e_1^1\cdot f_1^0=\mu f_1^1$ & null \\
\cline{3-4}
& & $e_1^1\cdot f_1^1=\mu f_2^0$ & $\rho(f_1^0)(e_1^1)=e_1^1$ \\
\cline{2-4}
& & trivial & $\rho(f_1^0)(e_1^1)=\lambda e_1^1,~\rho(f_2^0)(e_1^1)=\mu e_1^1$ \\
\cline{3-4}
& $L_{2|1}^6$ & $e_1^1\cdot f_1^1=\lambda f_1^0+\mu f_2^0$ & null \\
\cline{3-4}
& & $e_1^1\cdot f_1^0=\lambda f_1^1,~e_1^1\cdot f_2^0=\mu f_1^1$ & null \\
\hline
\end{tabular}
\caption{$(1|1, 2|1)$-type, $(\lambda, \mu\in\C)$}
\end{table}

\begin{table}
\resizebox{\textwidth}{!}{%
\parbox{1.35\textwidth}{%
\centering
\begin{tabular}{|c|c|l|l|}
\hline
$\textbf{A}$ & $\textbf{L}$ & $\textbf{Action}$ & $\textbf{Anchor}$\\
\hline
\multirow{35}{*}{$A_{1|1}^1$} & \multirow{2}{*}{$L_{2|2}^1$} & trivial & $\rho(f_1^0)(e_1^1)=\lambda e_1^1,~\rho(f_2^0)(e_1^1)=\mu e_1^1$ \\
\cline{3-4}
& & $e_1^1\cdot f_1^0=\lambda f_1^1,~e_1^1\cdot f_2^0=\mu f_1^1$ & $\rho(f_1^0)(e_1^1)= e_1^1$ \\
\cline{2-4}
& & trivial & $\rho(f_1^0)(e_1^1)=\lambda e_1^1,~\rho(f_2^0)(e_1^1)=\mu e_1^1$ \\
\cline{3-4}
& & $e_1^1\cdot f_1^0=\lambda(-if_1^1+f_2^1)$ & $\rho(f_1^0)(e_1^1)= e_1^1,~\rho(f_2^0)(e_1^1)=i e_1^1$ \\
\cline{3-4}
& $L_{2|2}^2$ & $e_1^1\cdot f_1^0=\lambda(if_1^1+f_2^1),~e_1^1\cdot f_2^0=\mu(if_1^1+f_2^1)$ & $\rho(f_1^0)(e_1^1)= e_1^1,~\rho(f_2^0)(e_1^1)=-i e_1^1$\\
\cline{3-4}	
& & $e_1^1\cdot f_1^0=\lambda(f_1^1-if_2^1),~e_1^1\cdot f_2^0=\mu(f_1^1-if_2^1)$ & $\rho(f_1^0)(e_1^1)= e_1^1,~\rho(f_2^0)(e_1^1)=-i e_1^1$ \\
\cline{3-4}
& & $e_1^1\cdot f_1^0=\lambda(f_1^1+if_2^1),~e_1^1\cdot f_2^0=\mu(f_1^1+if_2^1)$ & $\rho(f_1^0)(e_1^1)= e_1^1,~\rho(f_2^0)(e_1^1)=i e_1^1$ \\
\cline{2-4}
& \multirow{3}{*}{$L_{2|2}^3$} & trivial & $\rho(f_1^0)(e_1^1)= \lambda e_1^1$ \\
\cline{3-4}
& & $e_1^1\cdot f_1^0=\lambda f_2^1$ & $\rho(f_1^0)(e_1^1)=qe_1^1$ \\
\cline{3-4}
& & $e_1^1\cdot f_1^0=\lambda f_1^1$ & $\rho(f_1^0)(e_1^1)=pe_1^1$ \\
\cline{2-4}
& \multirow{2}{*}{$L_{2|2}^4$} & trivial & $\rho(f_1^0)(e_1^1)= \lambda e_1^1$ \\
\cline{3-4}
& & $e_1^1\cdot f_1^0=\lambda f_1^1$ & $\rho(f_1^0)(e_1^1)=pe_1^1$ \\
\cline{2-4}
& & trivial & $\rho(f_1^0)(e_1^1)=\lambda e_1^1$ \\
\cline{3-4}
& & $e_1^1\cdot f_1^0=\lambda(f_1^1-if_2^1)$& $\rho(f_1^0)(e_1^1)=(p-iq) e_1^1$ \\
\cline{3-4}
& $L_{2|2}^5$ & $e_1^1\cdot f_1^0=\lambda(f_1^1+if_2^1)$& $\rho(f_1^0)(e_1^1)=(p+iq) e_{1}^1$ \\
\cline{3-4}	
& & $e_1^1\cdot f_1^0=\lambda(if_1^1+f_2^1)$& $\rho(f_1^0)(e_1^1)=(p-iq) e_1^1$ \\
\cline{3-4}
& & $e_1^1\cdot f_1^0=\lambda(-if_1^1+f_2^1)$& $\rho(f_1^0)(e_1^1)=(p+iq) e_1^1$ \\
\cline{2-4}
& & trivial & $\rho(f_1^0)(e_1^1)=\lambda e_1^1$ \\
\cline{3-4}
& $L_{2|2}^6$ & $e_1^1\cdot f_1^0=\lambda f_2^1,~e_1^1\cdot f_2^0=-\lambda f_1^1$ & $\rho(f_1^0)(e_1^1)=pe_1^1$ \\
\cline{3-4}
& & $e_1^1\cdot f_1^0=\lambda f_1^1$ & $\rho(f_1^0)(e_1^1)=(1+p) e_1^1$ \\
\cline{2-4}
& & trivial & $\rho(f_1^0)(e_1^1)=\lambda e_1^1$ \\
\cline{3-4}
& $L_{2|2}^7$ &$e_1^1\cdot f_1^0=\lambda(f_1^1-if_2^1),~e_1^1\cdot f_1^1=2\lambda f_2^0$ &\multirow{2}{*}{$\rho(f_1^0)(e_1^1)=\frac{1}{2} e_1^1$} \\&&$e_1^1\cdot f_2^1=\pm2i\lambda f_2^0$&\\
\cline{2-4}
& \multirow{2}{*}{$L_{2|2}^8$} & trivial & $\rho(f_1^0)(e_1^1)=\lambda e_1^1$ \\
\cline{3-4}
& & $e_1^1\cdot f_1^0=\mu f_2^1$ & $\rho(f_1^0)(e_1^1)=\frac{1}{2} e_{1}^1$ \\
\cline{2-4}
& & trivial & $\rho(f_1^0)(e_1^1)=\lambda e_1^1$ \\
\cline{3-4}
& $L_{2|2}^9$ & $e_1^1\cdot f_1^0=\lambda f_2^1,~e_1^1\cdot f_1^1=\frac{1}{p} f_2^0$ & $\rho(f_1^0)(e_1^1)=(1-p)e_1^1$ \\
\cline{3-4}
& & $e_1^1\cdot f_1^0=(1-p)\lambda f_1^1,~e_1^1\cdot f_2^1=\lambda f_2^0$ & $\rho(f_1^0)(e_1^1)=pe_1^1$ \\
\cline{2-4}	
\multirow{9}{*}{$A_{1|1}^{1}$}	& \multirow{2}{*}{$L_{2|2}^{10}$} & trivial & $\rho(f_1^0)(e_1^1)=\lambda e_1^1$ \\
\cline{3-4}
& & $e_1^1\cdot f_1^0=\mu f_1^1$ & $\rho(f_1^0)(e_1^1)=\frac{1}{2} e_1^1$ \\
\cline{2-4}
& & trivial & $\rho(f_1^0)(e_1^1)=\lambda e_1^1$ \\
\cline{3-4}
& $L_{2|2}^{11}$ & $e_1^1\cdot f_1^0=\lambda(-if_1^1+f_2^1),~e_1^1\cdot f_1^1=\frac{2\lambda}{i+2p}f_2^0$ & \multirow{2}{*}{$\rho(f_1^0)(e_1^1)=\frac{1+2ip}{2} e_1^1$} \\&& $e_1^1\cdot f_2^1=\frac{2i\lambda}{i+2p}f_2^0$&\\
\cline{3-4}
& \multirow{2}{*}{$L_{2|2}^{11}$} & $e_1^1\cdot f_1^1=\frac{2\lambda}{2p-i}f_2^0,~e_1^1\cdot f_1^0=\lambda(if_1^1+f_2^1)$ & \multirow{2}{*}{$\rho(f_1^0)(e_1^1)=-\frac{2ip-1}{2} e_1^1$} \\&& $e_1^1\cdot f_2^1=\frac{2i\lambda}{2p-i}f_2^0$&\\
\cline{2-4}
& $L_{2|2}^{13}$ & $e_1^1\cdot f_1^1=\lambda f_1^0+\mu f_2^0,~e_1^1\cdot f_2^1=\gamma f_1^0+\theta f_2^0$ & null \\
\cline{2-4}
& \multirow{4}{*}{$L_{2|2}^{14}$} & trivial & $\rho(f_1^0)(e_1^1)=\lambda e_1^1$ \\
\cline{3-4}
& & $e_1^1\cdot f_2^1=\lambda f_2^0$ & $\rho(f_1^0)(e_1^1)=e_1^1$ \\
\cline{3-4}
& & $e_1^1\cdot f_1^0=\lambda f_2^1$ & \multirow{2}{*}{null} \\&&$e_1^1\cdot f_1^1=\lambda f_2^0$&\\
\cline{2-4}
& & trivial & $\rho(f_1^0)(e_1^1)=e_1^1$ \\
\cline{3-4}
& $L_{2|2}^{15}$ & $e_1^1\cdot f_2^0=\lambda f_2^1$ & $\rho(f_1^0)(e_1^1)=e_1^1$ \\
\cline{3-4}
& & $e_1^1\cdot f_2^1=\lambda f_2^0$ & null \\
\cline{2-4}
& & trivial & $\rho(f_1^0)(e_1^1)=\lambda e_1^1$ \\
\cline{3-4}
& & $e_1^1\cdot f_1^0=\lambda f_2^1,~e_1^1\cdot f_2^0=\mu f_2^1$ &\multirow{2}{*}{$\rho(f_1^0)(e_1^1)=- e_1^1,~\rho(f_1^1)(e_1^1)=\mu e_1^0$} \\&$L_{2|2}^{16}$&$e_1^1\cdot f_1^1=\lambda f_2^0-\mu f_1^0$&\\
\cline{3-4}
& & $e_1^1\cdot f_1^0=\lambda f_1^1,~e_1^1\cdot f_2^0=\mu f_1^1$ &\multirow{2}{*}{$\rho(f_1^0)(e_1^1)= e_1^1,~\rho(f_1^1)(e_1^1)=\mu e_1^0$} \\&&$e_1^1\cdot f_1^1=\mu f_1^0-\lambda f_2^0$&\\
\cline{2-4}
& & trivial & $\rho(f_1^0)(e_1^1)=\lambda e_1^1$ \\
\cline{3-4}
& $L_{2|2}^{17}$ & $e_1^1\cdot f_1^0=\lambda f_1^1,~e_1^1\cdot f_2^0=\mu f_1^1$ & \multirow{2}{*}{null} \\&&$e_1^1\cdot f_2^1=\mu f_1^0-\lambda f_2^0$&\\
\cline{2-4}
& \multirow{2}{*}{$L_{2|2}^{18}$} & trivial & $\rho(f_1^0)(e_1^1)= \lambda e_1^1,~\rho(f_2^0)(e_1^1)=\mu e_1^1$\\
\cline{3-4}
&& too many compatible actions to be listed & null\\
\hline							
\end{tabular}}}
\caption{$(1|1,2|2)$-type, $(\lambda, \mu, \gamma,\theta\in\C)$}
\end{table}

\begin{table}
\centering
\begin{tabular}{|c|c|l|l|}
\hline
$\textbf{A}$ & $\textbf{L}$ & $\textbf{Action}$ & $\textbf{Anchor}$ \\
\hline
\multirow{6}{*}{$A_{1|1}^1$} & \multirow{2}{*}{$L_{3|0}^1$} &\multirow{2}{*}{trivial} & $\rho(f_1^0)(e_1^1)=\lambda e_1^1,~\rho(f_2^0)(e_1^1)=\mu e_1^1$ \\&&&$\rho(f_3^0)(e_1^1)=\gamma e_1^1$\\
\cline{2-4}
& $L_{3|0}^2$ & trivial & $\rho(f_1^0)(e_1^1)=\lambda e_1^1,~\rho(f_2^0)(e_1^1)=\mu e_1^1$ \\
\cline{2-4}
& $L_{3|0}^3$ & trivial & $\rho(f_2^0)(e_1^1)=\lambda e_1^1,~\rho(f_3^0)(e_1^1)=\mu e_1^1$ \\
\cline{2-4}
& $L_{3|0}^4$ & trivial & $\rho(f_1^0)(e_1^1)=\lambda e_1^1$ \\
\cline{2-4}
& $L_{3|0}^5$ & trivial & $\rho(f_1^0)(e_1^1)=\lambda e_1^1$ \\
\hline	
\end{tabular}
\caption{$(1|1, 3|0)$-type, $(\lambda, \mu\in\C)$}
\end{table}

\begin{table}
\centering
\begin{tabular}{|c|c|l|l|}
\hline
$\textbf{A}$ & $\textbf{L}$ & $\textbf{Action}$ & $\textbf{Anchor}$ \\
\hline
& & trivial & $\rho(f_2^0)(e_1^1)=\lambda e_1^1,~\rho(f_3^0)(e_1^1)=\mu e_1^1$ \\
\cline{3-4}
\multirow{12}{*}{$A_{1|1}^1$}& $L_{3|1}^1$ & $e_1^1\cdot f_2^0=\lambda f_1^1$ & \multirow{2}{*}{$\rho(f_2^0)(e_1^1)= e_1^1$} \\&&$e_1^1\cdot f_3^0=\mu f_1^1$&\\
\cline{2-4}
&\multirow{2}{*}{$L_{3|1}^2$}&trivial&$\rho(f_3^0)(e_1^1)=\mu e_1^1$\\
\cline{3-4}
&& $e_1^1\cdot f_3^0=\lambda f_1^1$ &$\rho(f_3^0)(e_1^1)=qe_1^1$\\
\cline{2-4}
&\multirow{2}{*}{$L_{3|1}^3$}&trivial&$\rho(f_3^0)(e_1^1)=\mu e_1^1$\\
\cline{3-4}
&& $e_1^1\cdot f_3^0=\lambda f_1^1$ &$\rho(f_3^0)(e_1^1)=qe_1^1$\\
\cline{2-4}
& \multirow{2}{*}{$L_{3|1}^4$}& \multirow{1}{*}{trivial} & $\rho(f_2^0)(e_1^1)=\lambda e_1^1,~\rho(f_3^0)(e_1^1)=\mu e_1^1$ \\
\cline{3-4}
&& $e_1^1\cdot f_1^1=\lambda f_1^0$ & null \\
\cline{2-4}
&$L_{3|1}^5$& trivial & $\rho(f_1^0)(e_1^1)=\mu e_1^1$ \\
\cline{2-4}
&$L_{3|1}^6$& trivial & $\rho(f_1^0)(e_1^1)=\mu e_1^1$ \\
\cline{2-4}
& \multirow{5}{*}{$L_{3|1}^7$}& \multirow{2}{*}{trivial} & $\rho(f_1^0)(e_1^1)=\lambda e_1^1,~\rho(f_2^0)(e_1^1)=\mu e_1^1$ \\&&&$\rho(f_3^0)(e_1^1)=\gamma e_1^1$\\
\cline{3-4}
&&$e_1^1\cdot f_1^1=\lambda f_1^0+\mu f_2^0+\gamma f_3^0$&null\\
\cline{3-4}
&&$e_1^1\cdot f_1^1=\lambda f_1^0,~e_1^1\cdot f_2^0=\mu f_2^0$& \multirow{2}{*}{null}\\&&$e_1^1\cdot f_3^0=\gamma f_1^1$&\\
\hline
\end{tabular}
\caption{$(1|1, 3|1)$-type, $(\lambda, \mu, \gamma\in\C)$}
\end{table}

\begin{table}
\centering
\begin{tabular}{|c|c|l|l|}
\hline
$\textbf{A}$ & $\textbf{L}$ & $\textbf{Action}$ & $\textbf{Anchor}$ \\
\hline
\multirow{14}{*}{$A_{1|1}^1$} & \multirow{2}{*}{$L_{4|0}^1$} & \multirow{2}{*}{trivial} & $\rho(f_1^0)(e_1^1)=\lambda e_1^1,~\rho(f_2^0)(e_1^1)=\mu e_1^1$ \\&&&$\rho(f_3^0)(e_1^1)=\gamma e_1^1,~\rho(f_4^0)(e_1^1)=\theta e_1^1$\\
\cline{2-4}
& \multirow{2}{*}{$L_{4|0}^2$} & \multirow{2}{*}{trivial} & $\rho(f_1^0)(e_1^1)=\lambda e_1^1,~\rho(f_2^0)(e_1^1)=\mu e_1^1$ \\&&&$\rho(f_4^0)(e_1^1)=\gamma e_1^1$\\
\cline{2-4}
& \multirow{2}{*}{$L_{4|0}^3$} & \multirow{2}{*}{trivial} & $\rho(f_2^0)(e_1^1)=\lambda e_1^1,~\rho(f_3^0)(e_1^1)=\mu e_1^1$ \\&&&$\rho(f_4^0)(e_1^1)=\gamma e_1^1$\\
\cline{2-4}
& $L_{4|0}^4$ & trivial & $\rho(f_1^0)(e_1^1)=\lambda e_1^1,~\rho(f_4^0)(e_1^1)=\mu e_1^1$ \\
\cline{2-4}
& \multirow{2}{*}{$L_{4|0}^5$} & \multirow{2}{*}{trivial} & $\rho(f_1^0)(e_1^1)=\lambda e_1^1,~\rho(f_3^0)(e_1^1)=\mu e_1^1$ \\&&&$\rho(f_4^0)(e_1^1)=\gamma e_1^1$\\
\cline{2-4}
& $L_{4|0}^6$ & trivial & $\rho(f_2^0)(e_1^1)=\lambda e_1^1,~\rho(f_4^0)(e_1^1)=\mu e_1^1$ \\	
\cline{2-4}
& $L_{4|0}^7$ & trivial & $\rho(f_4^0)(e_1^1)=\lambda e_1^1$ \\	
\cline{2-4}		
& $L_{4|0}^8$ & trivial & $\rho(f_1^0)(e_1^1)=\lambda e_1^1,~\rho(f_2^0)(e_1^1)=\mu e_1^1$ \\
\cline{2-4}
& $L_{4|0}^k$ & \multirow{2}{*}{trivial} & \multirow{2}{*}{$\rho(f_1^0)(e_1^1)=\lambda e_1^1$} \\&$9\leq k\leq 16$&&\\	
\hline
\end{tabular}
\caption{$(1|1, 4|0)$-type, $(\lambda, \mu, \gamma,\theta\in\C)$}
\end{table}

$(2|0, 0|3)$ and $(2|0, 0|4)$-type: We already know that the anchor vanishes. There are too many suitable actions to be listed here (90 just for the pair $(\mathbf{A_{2|0}^2}, \mathbf{L_{0|3}^1}),$ for example).

\begin{table}
\centering
\begin{tabular}{|c|c|l|c|}
\hline
$\textbf{A}$ & $\textbf{L}$ & $\textbf{Action}$ & $\textbf{Anchor}$ \\
\hline
\multirow{9}{*}{$A_{2|0}^1$} & $L_{1|2}^1$ & trivial & $\rho(f_1^0)(e_2^0)=\lambda e_2^0$ \\
\cline{2-4}
& \multirow{2}{*}{$L_{1|2}^2$} & trivial& $\rho(f_1^0)(e_2^0)=\lambda e_2^0$ \\
\cline{3-4}
&&$e_2^0\cdot f_2^1=\lambda f_1^1$&null\\
\cline{2-4}
& \multirow{1}{*}{$L_{1|2}^3~\&~L_{1|2}^4$} & trivial& $\rho(f_1^0)(e_2^0)=\lambda e_2^0$ \\
\cline{2-4}
& \multirow{4}{*}{$L_{1|2}^6$} & trivial& $\rho(f_1^0)(e_2^0)=\lambda e_2^0$ \\
\cline{3-4}
&&$e_2^0\cdot f_2^1=\lambda f_1^1$& null\\
\cline{3-4}
&&$e_2^0\cdot f_1^1=-\lambda f_1^1+\mu f_2^1,~\mu\neq0$& \multirow{2}{*}{null}\\&&$e_2^0\cdot f_2^1=-\frac{\lambda^2}{\mu} f_1^1+\mu f_2^1$&\\
\hline
\multirow{10}{*}{$A_{2|0}^2$} & $L_{1|2}^k$ & $e_2^0\cdot f_1^0=f_1^0,~e_2^0\cdot f_1^1=f_1^1,~$ & \multirow{2}{*}{null} \\&$1\leq k\leq 5$&$e_2^0\cdot f_2^1=f_2^1$&\\
\cline{2-4}
&\multirow{8}{*}{$L_{1|2}^6$}& $e_2^0\cdot f_1^1=(1-\lambda)f_1^1+\mu f_2^1,~\mu\neq 0$&\multirow{8}{*}{null}\\&&$e_2^0\cdot f_2^1=\frac{\lambda(\lambda-1)}{\mu}f_1^1+\mu f_2^1$&\\
\cline{3-3}
&&$e_2^0\cdot f_2^1=\lambda f_1^1+f_2^1$&\\
\cline{3-3}
&&$e_2^0\cdot f_1^1=f_1^1,~e_2^0\cdot f_2^1=\lambda f_1^1$&\\
\cline{3-3}
&&$e_2^0\cdot f_1^0=f_1^0,~e_2^0\cdot f_2^1= f_2^1$&\\
\cline{3-3}
&&$e_2^0\cdot f_1^1=f_1^1,~e_2^0\cdot f_2^1=\lambda f_1^1$&\\
\cline{3-3}
&&$e_2^0\cdot f_1^0=f_1^0,~e_2^0\cdot f_1^1=f_1^1,~$&\\&&$e_2^0\cdot f_2^1=f_2^1$&\\
\hline
\end{tabular}
\caption{$(2|0, 1|2)$-type, $(\lambda,\mu \in\C)$}
\end{table}

\begin{table}
\centering
\begin{tabular}{|c|c|l|l|}
\hline
$\textbf{A}$ & $\textbf{L}$ & $\textbf{Action}$ & $\textbf{Anchor}$ \\
\hline
\multirow{9}{*}{$A_{2|0}^1$} & $L_{1|3}^k$ & \multirow{2}{*}{trivial} & \multirow{2}{*}{$\rho(f_1^0)(e_2^0)=\lambda e_2^0$} \\&$k\in \left\lbrace 1,3,4,6\right\rbrace $&&\\
\cline{2-4}
&\multirow{3}{*}{$L_{1|3}^2$}&trivial&$\rho(f_1^0)(e_2^0)=\lambda e_2^0$\\
\cline{3-4}
&&$e_2^0\cdot f_3^1=\lambda f_2^1$&null\\
\cline{3-4}
&&$e_2^0\cdot f_3^1=\lambda f_1^1$&$\rho(f_1^0)(e_2^0)=e_2^0$\\
\cline{2-4}
&\multirow{2}{*}{$L_{1|3}^5$}&trivial&$\rho(f_1^0)(e_2^0)=\lambda e_2^0$\\
\cline{3-4}
&&$e_2^0\cdot f_3^1=\lambda f_1^1$& null\\
\cline{2-4}
&\multirow{2}{*}{$L_{1|3}^8$}&trivial&$\rho(f_1^0)(e_2^0)=\lambda e_2^0$\\
\cline{3-4}
&& too many compatible actions to be listed&null\\
\hline
\multirow{3}{*}{$A_{2|0}^2$} & $L_{1|3}^k$ & $ e_2^0\cdot f_1^0=f_1^0,~e_2^0\cdot f_1^1=f_1^1$ & \multirow{2}{*}{null}\\&$1\leq k \leq 7$&$ e_2^0\cdot f_2^1=f_2^1,~e_2^0\cdot f_3^1=f_3^1$&\\
\cline{2-4}
&$L_{1|3}^8$& too many compatible actions to be listed & null\\
\hline
\end{tabular}
\caption{$(2|0, 1|3)$-type, $(\lambda\in\C)$}
\end{table}

\begin{table}
\resizebox{\textwidth}{!}{%
\parbox{1.1\textwidth}{%
\centering
\begin{tabular}{|c|c|l|c|}
\hline
$\textbf{A}$ & $\textbf{L}$ & $\textbf{Action}$ & $\textbf{Anchor}$ \\
\hline
\multirow{10}{*}{$A_{2|0}^1$} & \multirow{2}{*}{$L_{2|1}^1$} & trivial & $\rho(f_1^0)(e_2^0)=\lambda e_2^0$ \\
\cline{3-4}
&&$e_2^0\cdot f_1^0=\lambda f_2^0$&null\\
\cline{2-4}
& \multirow{2}{*}{$L_{2|1}^2$} & trivial & $\rho(f_1^0)(e_2^0)=\lambda e_2^0,~\rho(f_2^0)(e_2^0)=\lambda e_2^0$ \\
\cline{3-4}
&&$e_2^0\cdot f_1^0=\lambda (f_1^0-f_2^0),~e_2^0\cdot f_2^0=\lambda (f_1^0+f_2^0)$&null\\
\cline{2-4}
& \multirow{1}{*}{$L_{2|1}^k$} & trivial & $\rho(f_1^0)(e_2^0)=\lambda e_2^0$ \\
\cline{3-4}
&$3\leq k\leq 5$&$e_2^0\cdot f_1^0=\lambda f_2^0$&$\rho(f_1^0)(e_2^0)= e_2^0$\\
\cline{2-4}
& \multirow{4}{*}{$L_{2|1}^6$} & trivial & $\rho(f_1^0)(e_2^0)=\lambda e_2^0,~\rho(f_2^0)(e_2^0)=\lambda e_2^0$ \\
\cline{3-4}
&&$e_2^0\cdot f_1^0=\lambda f_2^0$&null\\
\cline{3-4}
&&$e_2^0\cdot f_1^0=-\lambda f_1^0+\mu f_2^0,~\mu\neq0$&\multirow{2}{*}{null}\\&&$e_2^0\cdot f_2^0=-\frac{\lambda^2}{\mu} f_1^0+\lambda f_2^0$&\\
\cline{1-4}
\multirow{23}{*}{$A_{2|0}^2$} & \multirow{4}{*}{$L_{2|1}^1$} &$e_2^0\cdot f_1^0=f_1^0+\lambda f_2^0$ & \multirow{4}{*}{null} \\
\cline{3-3}
& &$e_2^0\cdot f_1^0=\lambda f_2^0$ & \\
\cline{3-3}
& &$e_2^0\cdot f_2^0=f_2^0,~e_2^0\cdot f_1^1=f_1^1$ & \\
\cline{3-3}
& &$e_2^0\cdot f_1^0=f_1^0,~e_2^0\cdot f_2^0=f_2^0,~e_2^0\cdot f_1^1=f_1^1$ & \\
\cline{2-4}
& \multirow{6}{*}{$L_{2|1}^2$} &$e_2^0\cdot f_1^0=(1-\lambda)(f_1^0+ f_2^0)$ & \multirow{6}{*}{null} \\&&$e_2^0\cdot f_2^0=\lambda(f_1^0+f_2^0)$&\\
\cline{3-3}
& &$e_2^0\cdot f_1^0=f_1^0+f_2^0$ & \\
\cline{3-3}
& &$e_2^0\cdot f_1^0=(1-\lambda)f_1^0-\lambda f_2^0,$ & \\&&$e_2^0\cdot f_2^0=-(1-\lambda)f_1^0+\lambda f_2^0$&\\
\cline{3-3}
& &$e_2^0\cdot f_1^0=f_1^0,~e_2^0\cdot f_2^0=\pm f_2^0,~e_2^0\cdot f_1^1=f_1^1$ & \\
\cline{2-4}
& $L_{2|1}^3~ \& ~L_{2|1}^4$ &$e_2^0\cdot f_1^0=f_1^0,~e_2^0\cdot f_2^0=f_2^0,~e_2^0\cdot f_1^1=f_1^1$ & null \\
\cline{2-4}
& \multirow{3}{*}{$L_{2|1}^5$} &$e_2^0\cdot f_1^1=f_1^1$ & \multirow{3}{*}{null}\\
\cline{3-3}
&& $e_2^0\cdot f_1^0=f_1^0,~e_2^0\cdot f_2^0=f_2^0$ & \\
\cline{3-3}
& &$e_2^0\cdot f_1^0=f_1^0,~e_2^0\cdot f_2^0=f_2^0,~e_2^0\cdot f_1^1=f_1^1$ & \\
\cline{2-4}
& \multirow{1}{*}{$L_{2|1}^6$} &$e_2^0\cdot f_1^0=(1-\lambda)f_1^0+\mu f_2^0,~\mu\neq 0$ & \multirow{1}{*}{null} \\&&$e_2^0\cdot f_2^0=\frac{\lambda(\lambda-1)}{\mu}f_1^0+\lambda f_2^0$&\\
\cline{3-3}
&\multirow{5}{*}{$L_{2|1}^6$}&$e_2^0\cdot f_1^0=f_1^0+\lambda f_2^0$&\multirow{5}{*}{null}\\
\cline{3-3}
&&$e_2^0\cdot f_1^0=f_1^0+\lambda f_2^0,~e_2^0\cdot f_1^1=f_1^1$&\\
\cline{3-3}
&&$e_2^0\cdot f_1^0=f_1^0+\lambda f_2^0$&\\
\cline{3-3}
&&$e_2^0\cdot f_2^0=\lambda f_1^0+f_2^0,~e_2^0\cdot f_1^1=f_1^1$&\\
\cline{3-3}
&&$e_2^0\cdot f_2^0=\lambda f_1^0$&\\
\cline{3-3}
&&$e_2^0\cdot f_1^1=f_1^1$&\\
\cline{3-3}
&&$e_2^0\cdot f_1^0=f_1^0,~e_2^0\cdot f_2^0=f_2^0,~e_2^0\cdot f_1^1=f_1^1$&\\
\hline
\end{tabular}}}
\caption{$(2|0,2|1)$-type, $(\lambda,\mu\in\C)$}
\end{table}

\begin{table}
\resizebox{\textwidth}{!}{%
\parbox{1.2\textwidth}{%
\centering
\begin{tabular}{|c|c|l|l|}
\hline
$\textbf{A}$ & $\textbf{L}$ & $\textbf{Action}$ & $\textbf{Anchor}$ \\
\hline
\multirow{19}{*}{$A_{2|0}^1$} & \multirow{2}{*}{$L_{2|2}^1$} & trivial & $\rho(f_1^0)(e_2^0)=\lambda e_2^0,~\rho(f_2^0)(e_2^0)=\lambda e_2^0$ \\
\cline{3-4}
&&$e_2^0 \cdot f_1^0=\lambda f_2^0,~e_2^0 \cdot f_2^1=\lambda f_1^1 $ & null  \\
\cline{2-4}
& \multirow{1}{*}{$L_{2|2}^2$} & trivial & $\rho(f_1^0)(e_2^0)=\lambda e_2^0,~\rho(f_2^0)(e_2^0)=\lambda e_2^0$ \\
\cline{2-4}
& $L_{2|2}^k$ & trivial & $\rho(f_1^0)(e_2^0)=\lambda e_2^0$ \\
\cline{3-4}
&$3\leq k\leq 15,$ &\multirow{2}{*}{$e_2^0 \cdot f_1^0=\lambda f_2^0$} & \multirow{2}{*}{$\rho(f_1^0)(e_2^0)=e_2^0$}  \\&$k\notin \left\lbrace 6,12,13\right\rbrace$&&\\
\cline{2-4}
& \multirow{2}{*}{$L_{2|2}^6$} & trivial & $\rho(f_1^0)(e_2^0)=\lambda e_2^0$ \\
\cline{3-4}
&&$e_2^0 \cdot f_1^0=\lambda f_2^0,~e_2^0 \cdot f_2^1=\frac{\lambda}{p}f_1^1$& $\rho(f_1^0) (e_2^0)=e_2^0$\\
\cline{2-4}
& \multirow{2}{*}{$L_{2|2}^{16}$} & trivial & $\rho(f_1^0)(e_2^0)=\lambda e_2^0$ \\
\cline{3-4}
&&$e_2^0 \cdot f_1^0=\lambda f_2^0,$& null\\
\cline{2-4}
& \multirow{2}{*}{$L_{2|2}^{17}$} & trivial & $\rho(f_1^0)(e_2^0)=\lambda e_2^0$ \\
\cline{3-4}
&&$e_2^0 \cdot f_1^0=\lambda f_2^0,~e_2^0 \cdot f_2^1=\mu f_1^1$& null\\
\cline{2-4}
& \multirow{2}{*}{$L_{2|2}^{18}$} & trivial & $\rho(f_1^0)(e_2^0)=\lambda e_2^0,~\rho(f_2^0)(e_2^0)=\mu e_2^0$ \\
\cline{3-4}
&&$e_2^0 \cdot f_1^0=\lambda f_2^0,~e_2^0 \cdot f_2^1=\mu f_1^1$&\multirow{6}{*}{null}\\
\cline{3-3}
&\multirow{2}{*}{$(\mu,\theta \neq 0)\dashrightarrow $}&$e_2^0 \cdot f_1^0=-\lambda f_2^0+\mu f_2^0,~e_2^0 \cdot f_2^0=-\frac{\lambda^2}{\mu}f_1^0+\lambda f_2^0$&\\&&$e_2^0 \cdot f_1^1=-\gamma f_1^1+\theta f_2^1,~e_2^0 \cdot f_2^1=-\frac{\gamma^2}{\theta}f_1^1+\gamma f_2^1$&\\
\cline{3-3}
&&$e_2^0 \cdot f_2^0=\lambda f_1^0,~e_2^0 \cdot f_2^1=\mu f_1^1$&\\
\cline{3-3}
&$(\mu\neq 0)\dashrightarrow$&$e_2^0 \cdot f_1^1=-\lambda f_1^1+\mu f_2^1,~e_2^0 \cdot f_2^1=-\frac{\lambda^2}{\mu}f_1^1+\mu f_2^1$&\\
\cline{3-3}
&$(\mu\neq 0)\dashrightarrow$&$e_2^0 \cdot f_1^0=-\lambda f_1^0+\mu f_2^0,~e_2^0 \cdot f_2^0=-\frac{\lambda^2}{\mu}f_1^0+\mu f_2^0$&\\
\hline
\multirow{4}{*}{$A_{2|0}^2$} & \multirow{1}{*}{$L_{2|2}^k$} & \multirow{2}{*}{$e_2^0\cdot f_1^0=f_1^0,~e_2^0\cdot f_2^0=f_2^0,$} & \multirow{3}{*}{null} \\ &$1\leq k\leq 17$&\multirow{2}{*}{$e_2^0\cdot f_1^1=f_1^1,~e_2^0\cdot f_2^1=f_2^1$}&\\&$k\notin \left\lbrace 2,13,15\right\rbrace$ &&\\
\cline{2-4}				
\multirow{15}{*}{$A_{2|0}^2$}&\multirow{4}{*}{$L_{2|2}^2$}&$e_2^0 \cdot f_1^0=\frac{1}{2}f_1^0-\frac{i}{2}f_2^0,~e_2^0 \cdot f_2^0=\frac{i}{2}f_1^0+\frac{1}{2}f_2^0$&\multirow{4}{*}{null}\\&&$e_2^0 \cdot f_1^1=\frac{1}{2}f_1^1+\frac{i}{2}f_2^1,~e_2^0 \cdot f_2^1=-\frac{i}{2}f_2^1+\frac{1}{2}f_2^1$&\\
\cline{3-3}
&&$e_2^0 \cdot f_1^0=\frac{1}{2}f_1^0+\frac{i}{2}f_2^0,~e_2^0 \cdot f_2^0=-\frac{i}{2}f_1^0+\frac{1}{2}f_2^0$&\\
&&$e_2^0 \cdot f_1^1=\frac{1}{2}f_1^1-\frac{i}{2}f_2^1,~e_2^0 \cdot f_2^1=\frac{i}{2}f_2^1+\frac{1}{2}f_2^1$&\\
\cline{2-4}
&\multirow{2}{*}{$L_{2|2}^{13}$}&$e_2^0\cdot f_1^0=f_1^0,~e_2^0\cdot f_2^0=f_2^0,$&\multirow{4}{*}{null}\\&&$e_2^0\cdot f_1^1=f_1^1,~e_2^0\cdot f_2^1=f_2^1$&\\
\cline{3-3}
&&$e_2^0\cdot f_2^0=f_2^0,~e_2^0\cdot f_2^1=f_2^1$&\\
\cline{3-3}
&&$e_2^0\cdot f_1^0=f_1^0,~e_2^0\cdot f_1^1=f_1^1$&\\
\cline{2-4}
&\multirow{5}{*}{$L_{2|2}^{15}$}&$e_2^0\cdot f_1^0=f_1^0,~e_2^0\cdot f_2^0=f_2^0,$&\multirow{5}{*}{null}\\&&$e_2^0\cdot f_1^1=f_1^1,~e_2^0\cdot f_2^1=f_2^1$&\\
\cline{3-3}
&&$e_2^0\cdot f_2^1=f_2^1$&\\
\cline{3-3}
&&$e_2^0\cdot f_1^0=f_1^0,~e_2^0\cdot f_2^0=f_2^0$&\\&&$e_2^0\cdot f_1^1=f_1^1$&\\
\cline{2-4}
&\multirow{1}{*}{$L_{2|2}^{18}$}&too many compatible actions to be listed&null\\
\hline
\end{tabular}}}
\caption{$(2|0,2|2)$-type, $(\lambda, \mu, \gamma, \theta\in\C)$}
\end{table}

$(2|0, 3|0)$-type, $(2|0, 3|1)$-type and $(2|0, 4|0)$-type: no results to mention here.

\EditInfo{November 15, 2021}{November 25, 2022}{Friedrich Wagemann}

\end{paper}